\input amstex\documentstyle{amsppt}  
\pagewidth{12.5cm}\pageheight{19cm}\magnification\magstep1
\topmatter
\title From conjugacy classes in the Weyl group to unipotent classes, III\endtitle
\rightheadtext{Conjugacy classes in the Weyl group and unipotent classes, III}
\author G. Lusztig\endauthor
\address{Department of Mathematics, M.I.T., Cambridge, MA 02139}\endaddress
\thanks{Supported in part by the National Science Foundation}\endthanks
\endtopmatter   
\document
\define\Irr{\text{\rm Irr}}

\define\us{\un s}

\define\ul{\un l}

\define\uD{\un D}
\define\uG{\un G}

\define\uWW{\un\WW}

\define\mpb{\medpagebreak}

\define\dsv{\dashv}

\define\pe{\perp}
\define\si{\sim}

\define\sqc{\sqcup}

\define\qua{\quad}

\define\bG{\bar G}

\define\lb{\linebreak}

\define\bin{\binom}
\define\op{\oplus}
   
\redefine\sp{\spadesuit}
\define\part{\partial}
\define\em{\emptyset}

\define\ra{\rangle}
\define\n{\notin}

\define\m{\mapsto}
\define\do{\dots}
\define\la{\langle}

\define\lra{\leftrightarrow}

\define\sub{\subset}    

\define\T{\times}
\define\ti{\tilde}
\define\nl{\newline}
\redefine\i{^{-1}}

\define\un{\underline}

\define\ot{\otimes}
\define\bbq{\bar{\QQ}_l}

\define\Ad{\text{\rm Ad}}
\define\Hom{\text{\rm Hom}}
\define\End{\text{\rm End}}

\define\tr{\text{\rm tr}}

\define\a{\alpha}

\define\g{\gamma}
\redefine\d{\delta}
\define\e{\epsilon}
\define\et{\eta}

\redefine\o{\omega}
\define\p{\pi}
\define\ph{\phi}
\define\ps{\psi}
\define\r{\rho}
\define\s{\sigma}
\redefine\t{\tau}

\define\k{\kappa}
\redefine\l{\lambda}
\define\z{\zeta}
\define\x{\xi}

\redefine\G{\Gamma}
\redefine\D{\Delta}

\define\Si{\Sigma}

\redefine\L{\Lambda}
\define\Ph{\Phi}
\define\Ps{\Psi}

\define\kk{\bold k}

\define\nn{\bold n}

\define\ww{\bold w}

\define\NN{\bold N}

\define\QQ{\bold Q}

\define\WW{\bold W}
\define\ZZ{\bold Z}
\define\XX{\bold X}

\define\cb{\Cal B}

\define\cf{\Cal F}
\define\cg{\Cal G}
\define\ch{\Cal H}

\define\cl{\Cal L}

\define\co{\Cal O}

\define\cs{\Cal S}
\define\ct{\Cal T}

\define\cz{\Cal Z}

\define\fe{\frak e}

\define\fB{\frak B}

\define\fS{\frak S}

\define\tg{\ti g}

\define\ts{\ti s}

\define\tz{\ti z}

\define\tA{\ti A}

\define\tQ{\ti Q}

\define\tV{\ti V}

\define\sha{\sharp}

\define\sps{\supset}

\define\che{\check}

\define\chg{\che{g}}

\define\CA{Ca}
\define\DL{DL}
\define\GP{GP}
\define\GKP{GKP}
\define\HE{H1}
\define\HEA{H2}
\define\HL{HL}
\define\LU{L1}
\define\OR{L2}
\define\CSIV{L3}
\define\CDI{L4}
\define\CDII{L5}
\define\WE{L6}
\define\WEII{L7}
\define\WEH{L8}
\define\XW{L9}
\define\MA{M1}
\define\MAA{M2}
\define\MAB{M3}
\define\LX{LX}
\define\SP{Sp}

\head Introduction\endhead
\subhead 0.1\endsubhead
Let $G$ be an affine algebraic group over an algebraically closed field of characteristic $p\ge0$ such that the 
identity component $G^0$ of $G$ is reductive. Let $\WW$ be the Weyl group of $G^0$. In \cite{\WE} we have defined
(assuming that $G=G^0$) a natural (surjective) map $\Ph$ from the set $\uWW$ of conjugacy classes in $\WW$ to the
set $\uG$ of unipotent conjugacy classes in $G$. In this paper we extend the definition of $\Ph$ to the case 
where $G$ is not necessarily connected, by replacing $\uG$ by the set of unipotent $G^0$-conjugacy classes in a 
fixed connected component $D$ of $G$ whose image in $G/G^0$ is unipotent and $\uWW$ by the set of "twisted 
conjugacy classes" (in a suitable sense depending on $D$) of $\WW$. The general case can be reduced to the 
already known case when $G=G^0$ and to four other cases: the case where $G$ has two components ($G^0=PGL_n$ and 
$D$) and conjugation by some element in $D$ takes a matrix to its transpose inverse with $p=2$ (see \S4,\S5); the
case where $G$ is an even full orthogonal group and $D\ne G^0$ with $p=2$ (see \S3); and two exceptional cases 
related to $E_6,p=2$ and $D_4,p=3$ (which can be treated using in part computer calculations, see \S2). As a 
biproduct of our analysis in \S4 we obtain a new description of certain varieties from \cite{\DL} corresponding 
to a unitary group over a finite field (see \S6). On the other hand in 5.11 we show that a map similar to $\Ph$ 
can be defined in a case (arising from an outer automorphism of $GL_n$) where $D$ does not contain unipotent 
elements.

I thank Gongqin Li for her help with programming in GAP3.

\subhead 0.2. Notation\endsubhead 
For any set of integers $M$ we set $M_{ev}=\{x\in M;x\text{ even}\}$, $M_{odd}=\{x\in M;x\text{ odd}\}$. For any 
collection of vectors $e_1,e_2,\do,e_t$ in a vector space we denote by $S(e_1,e_2,\do,e_t)$ the subspace spanned 
by $e_1,e_2,\do,e_t$. For any group $\G$ let $\cz_\G$ be the centre of $\G$; if $g\in\G$ and $\G'$ is a subset of
$\G$ let $Z_{\G'}(g)=\{g'\in\G';gg'=g'g\}$.

\subhead 0.3\endsubhead
{\it Errata to \cite{\WE}.} In 4.4 replace

"We can find integers $a_1,a_2,\do,a_t,b_1,b_2,\do,b_t$ (all $\ge1$) such that ..." by:

"We can find integers $a_1,a_2,\do,a_t,b_1,b_2,\do,b_t$ (all $\ge1$) such that $b_1+b_2+\do+b_t=\s$,"

In 4.4 replace the equation $Y=(a_1+a_2+\do+a_t)(b_1^2-b_1)/2+\do$ by the equation
$$\align&Y=(a_1+a_2+\do+a_t)(b_1^2-b_1)/2\\&+(a_1+a_2+\do+a_{t-1})((b_1+b_2)^2-b_1-b_2)-(b_1^2-b_1))/2+\do+\\&
a_1((b_+\do+b_t)^2-b_1-\do-b_t)-((b_1+\do+b_{t-1})^2-b_1-\do-b_{t-1}))/2.\endalign$$
Replace 1.6(a) by 
$$\align&1\m2\m\do\m p_1\m\nn\m\nn-1\m\do\m\nn-p_1+1\m1,\\&
         p_1+1\m p_1+2\m\do\m p_1+p_2\m\nn-p_1\m\nn-p_1-1\m\do \\&  \m\nn-p_1-p_2+1\m p_1+1,\\&\do\\&
p_{<\s}+1\m p_{<\s}+2\m\do\m p_{<\s}+p_\s\m\nn-p_{<\s}\\&\m\nn-p_{<\s}-1\m\do\m\nn-p_{<\s}-p_\s+1\m p_{<\s}+1,\\&
\text{ and, if }\k=1,\\&p_{n+1}\m p_{n+1}.\tag a\endalign$$

\head Contents\endhead
1. The main results; statements and preliminary reductions.

2. The exceptional cases.

3. Even full orthogonal groups.

4. Bilinear forms.

5. Almost unipotent bilinear forms.

6. Finite unitary groups.

7. Final remarks.

\head 1. The main results; statements and preliminary reductions\endhead
\subhead 1.1\endsubhead
Let $\kk$ be an algebraically closed field of characteristic $p\ge0$. In this paper we fix a group $G$ with a 
given normal subgroup $G^0$ and a given left/right $G^0$-coset denoted by $D$; we assume that $G^0$ is a 
connected reductive algebraic group over $\kk$ and that for any $g\in D$ the map $x\m gxg\i$ is an isomorphism of
algebraic groups $G^0@>\si>>G^0$. Let $\cz_{G^0}^D$ be the set of all $z\in\cz_{G^0}$ such that $zg=gz$ for 
some/any $g\in D$; this is a closed subgroup of $\cz_{G^0}$. Let $\cb$ be the variety of Borel subgroups 
of $G^0$. Let $\WW$ be the Weyl group of $G^0$. We view $\WW$ as a set indexing the set of orbits of $G^0$ acting
on $\cb\T\cb$ by $g:(B,B')\m(gBg\i,gB'g\i)$. For $w\in\WW$ we write $\co_w$ for the corresponding $G^0$-orbit. 
The group structure on $\WW$ is defined as in \cite{\WE, 0.2}. Define $\ul:\WW@>>>\NN$ by 
$\ul(w)=\dim\co_w-\dim\cb$ (length function). Let $S=\{s\in\WW;\ul(s)=1\}$. For any $J\sub S$ let $\WW_J$ be the 
subgroup of $\WW$ generated by $J$. Now $D$ defines a group automorphisms $\e_D:\WW@>>>\WW$ (preserving length) 
by the requirement that
$$(B,B')\in\co_w,g\in D\implies(gBg\i,gB'g\i)\in\co_{\e_D(w)}.$$
Let $\uWW_D$ be the set of $\e_D$-conjugacy classes in $\WW$ that is, the set of equivalence classes in $\WW$ for 
the orbits of the $\WW$-action $w_1:w\m w_1\i w\e_D(w_1)$ on $\WW$. An element $w\in\WW$ (or its $\e_D$-conjugacy
class $C$) is said to be $\e_D$-elliptic if for any $J\sub S$ such that $\e_D(J)=J$, $J\ne S$, we have 
$C\cap\WW_J=\em$. Let $\uWW_D^{el}$ be the set of elliptic $\e_D$-conjugacy classes in $\WW$. Let $\WW_{D-min}$ 
be the set of all $w\in\WW$ such that $\ul(w)\le\ul(w')$ for any $w'$ in the $\e_D$-conjugacy class of $w$. For 
any $\e_D$-conjugacy class $C$ in $\WW$ we set $C_{min}=C\cap\WW_{D-min}$.

In the remainder of this section we assume that 

(a) $G$ is an affine algebraic group over $\kk$ and that $G^0$ is its identity component (so that $D$ is a
connected component of $G$)
\nl
and that

(b) the image of $D$ in $G/G^0$ is unipotent (that is, its order is a power of $p$ if $p>1$ and $1$ if $p=0$).
\nl
Let $\uD$ be the set of $G^0$-conjugacy classes in $D$ which are unipotent. We say that $\g\in\uD$ is 
distinguished if for some/any $g\in\g$, any torus in $Z_G(g)$ is contained in $\cz_{G^0}^D$. Let 
$\uD_{dis}=\{g\in\uD;\g\text{ distinguished}\}$. For $w\in\WW,\g\in\uD$ we write $w\dsv_D\g$ if for some/any 
$g\in\g$ we have $\{B\in\cb;(B,gBg\i)\in\co_w\}\ne\em$. For $w\in\WW$ we set $\Si_{w,D}=\{\g\in\uD; w\dsv_D\g\}$;
we regard $\Si_{w,D}$ as a partially ordered set where $\g\le\g'$ if $\g\sub\bar\g'$ (closure of $\g'$ in $D$).

\subhead 1.2\endsubhead
We shall need the following result.

(a) {\it Assume that $w,w'$ are elements of $\WW_{D-min}$ which are $\e_D$-conjugate. Then 
$\Si_{w,D}=\Si_{w',D}$.}
\nl
To prove this we can assume that $\kk$ is an algebraic closure of a finite field $F_q$ with $q$ elements and that
$G$ has a fixed $F_q$-rational structure such that $G^0$ is $F_q$-split, $D$ is defined over $F_q$ and each 
unipotent $G^0$-orbit in $D$ is defined over $F_q$. Let $\g\in\uD$. Let $g\in\g(F_q)$. Let $N_q$ (resp. $N'_q$) 
be the number of $F_q$-rational points of the $F_q$-variety $\{B\in\cb;(B,gBg\i)\in\co_w\}$ (resp. 
$\{B\in\cb;(B,gBg\i)\in\co_{w'}\}$). Note that $\g\in\Si_{w,D}$ if and only if $N_{q^s}\ne0$ for some $s\ge1$; 
similarly, $\g\in\Si_{w',D}$ if and only if $N'_{q^s}\ne0$ for some $s\ge1$. It is then enough to show that 
$N_{q^s}=N'_{q^s}$ for all $s\ge1$. Replacing $q^s$ by $q$ we can also assume that $s=1$. Let $\cf$ be the set of
functions $\cb(F_q)@>>>\bbq$ ($l$ is a fixed prime number $\ne p$). For any $x\in G(F_q)$ we define 
$S_x:\cf@>>>\cf$ by $f\m f'$ where $f'(B)=f(x^{-1}Bx)$. For any $y\in\WW$ we denote by $T_y:\cf@>>>\cf$ the 
linear map $f\m f'$ where 
$$f'(B)=\sum_{B'\in\cb(F_q);(B,B')\in\co_y}f(B').$$
Let $\ch_q$ be the subspace of $\End(\cf)$ spanned by
$T_y(y\in\WW$); this is a subalgebra of $\End(\cf)$. From the definitions we have $N_q=\tr(T_wS_g:\cf@>>>\cf)$, 
$N'_q=\tr(T_{w'}S_g:\cf@>>>\cf)$. Thus it is enough to show that $\tr((T_w-T_{w'})S_g:\cf@>>>\cf)=0$. We define a
linear map $\mu:\ch_q@>>>\ch_q$ by $\mu(T_y)=T_{\e_D(y)}$ for all $y\in\WW$; this is an algebra automorphism. 
From the definitions for any $y\in\WW$ we have $T_{\e_D(y)}S_g=S_gT_y:\cf@>>>\cf$ (we use that $g\in D$); hence 
for any $\x\in\ch_q$ we have $\mu(\x)S_g=S_g\x$. From our assumptions on $w,w'$ and from \cite{\GKP, 7.2(a)} we 
see that $T_w-T_{w'}$ is a linear combinations of elements of the form $\x\x'-\x'\mu(\x)$ ($\x,\x'\in\ch_q$). 
Hence it is enough to show that $\tr((\x\x'-\x'\mu(\x))S_g:\cf@>>>\cf)=0$ for any $\x,\x'\in\ch_q$. The last 
trace is equal to $\tr(\x\x'S_g-\x'S_g\x:\cf@>>>\cf)$ and this is clearly $0$. This completes the proof of (a).

\mpb

In view of (a), for any $C\in\uWW_D$ we can define $\Si_{C,D}=\Si_{w,D}$ where $w$ is any element of $C_{min}$.

The following is one of the main results of this paper.

\proclaim{Theorem 1.3}(a) Let $C\in\uWW_D$. There exists (a necesarily unique) $\g\in\Si_{C,D}$ such that
$\g\le\g'$ for all $\g'\in\Si_{C,D}$. We set $\g=\Ph(C)$.

(b) If $C,C'$ are elements of $\uWW_D^{el}$ and $\Ph(C)=\Ph(C')$, then $C=C'$.

(c) For any $\g\in\uD$ there exists $C\in\uWW_D$ such that $\g=\Ph(C)$.

(d) For any $\g\in\uD_{dis}$ there exists $C\in\uWW_D^{el}$ such that $\g=\Ph(C)$.
\endproclaim
We now state a variant (in fact a special case) of the Theorem.
\proclaim{Proposition 1.4} 1.3(a) holds for any $C\in\uWW_D^{el}$; moreover, 1.3(b) and 1.3(d) hold.
\endproclaim
We now show:

(a) {\it If the proposition holds for $G,D$ and for any $G',D'$ such that $\dim(G')<\dim G$ then Theorem 1.3 
holds for $G,D$.}
\nl
The proof is a generalization of that in \cite{\WE, 1.1}. For any $J\sub S$ such that $\e_D(J)=J$ we denote by 
$P_J$ a parabolic subgroup of type $J$ of $G^0$, by $L_J$ a Levi subgroup of $P_J$ and by $U_{P_J}$ the unipotent
radical of $P_J$; we identify $\WW_J$ with the Weyl group of $L_J$ in the standard way (using $P_J$). By 
\cite{\CDI, 1.26} we have $N_GP_J=(N_GL_J\cap N_GP_J)U_{P_J}$, $(N_GL_J\cap N_GP_J)\cap U_{P_J}=\{1\}$. Here 
$N_G()$ denotes normalizer in $G$. It follows that $N_DP_J=(N_DL_J\cap N_DP_J)U_{P_J}$ where $N_D()=N_G()\cap D$.
Since $\e_D(J)=J$ we have $N_DP_J\ne\em$ hence $D_J:=N_DL_J\cap N_DP_J\ne\em$. Note that $D_J$ is a connected 
component of $N_GL_J$ (whose identity component is $L_J$) and $(\uWW_J)_{D_J}$, $\uD_J$ are defined in terms of 
$N_GL_J,D_J$ in the same way as $\uWW_D$, $\uD$ are defined in terms of $G,D$. Now let $C\in\uWW_D$. We can find 
$J\sub S$ as above and $C'\in(\uWW_J)_{D_J}^{el}$ such that $C'=C\cap\WW_J$. (For future reference we denote by 
$\mu(C)$ the number of $\e_D$-orbits on $S-J$; it is independent of the choice of $C$. Note that $C$ is 
$\e_D$-elliptic if and only if $\mu(C)=0$.) We set $P_J=P,L_J=L$, $U_{P_J}=U$. 
By our assumption, $\g_0:=\Ph(C')$ is a well defined unipotent $L$-conjugacy class in $D_J$. Let $\g$ be the 
unipotent $G^0$-conjugacy class in $D$ that contains $\g_0$. Let $w\in C'_{min}$. Let $\g'\in\uD$ be such that 
for some $B,B'\in\cb$, $g'\in\g'$ we have $(B,B')\in\co_w$, $B'=g'Bg'{}\i$. Replacing $B,B',g'$ by 
$xBx\i,xB'x\i,xg'x\i$ for some $x\in G^0$ we see that we can assume that $B\sub P$ and then we automatically have
$B'\sub P$ that is, $g'Bg'{}\i\sub P$. We have also $g'Bg'{}\i\sub g'Pg'{}\i$ and $g'Pg'{}\i$ is of the same type
as $P$ (since $\e_D(J)=J$) hence $g'Pg'{}\i=P$ that is, $g'\in N_DP$. We have $g'=g'_1v$ where $g'_1\in D_J$ is 
unipotent and $v\in U$. We shall use the following fact which will be verified below.

(b) {\it For any element $h\in D_J$ and any $v'\in U$ we can find a one parameter subgroup $\l:\kk^*@>>>\cz_L$ 
such that $\l(t)h=h\l(t)$ for all $t$ and $\l(t)v'\l(t\i)$ converges to $1$ when $t\in\kk^*$ converges to $0$.}
\nl
Using (b) with $h=g'_1,v'=v$ we see that $\l(t)g'\l(t)\i=g'_1\l(t)v\l(t)\i$ converges to $g'_1$ when $t\in\kk^*$ 
converges to $0$. Thus $g'_1$ is contained in the closure of $\g'$. Hence the $L$-conjugacy class of $g'_1$ is 
contained in the closure of $\g'$. Note also that $B'=g'_1Bg'_1{}\i$. By definition, $\g_0$ is contained in the 
closure of the $L$-conjugacy class of $g'_1$. Hence $\g_0$ is contained in the closure of $\g'$ and $\g$ is 
contained in the closure of $\g'$. We see that $\Si_{w,D}$ has a unique minimal element namely $\g$. Since 
$w\in C_{min}$, we see that 1.3(a) holds for $G,D$ (assuming (b)).

We now prove (b). Let $u$ be a unipotent, quasi-semisimple element in $D_J$. We have $h=ub$ where $b\in L$ hence
$b$ commutes with any element in $\cz_L$. Hence it is enough to prove (b) with $h$ replaced by $u$. We can find a 
Borel subgroup $B_1$ of $L$ and a maximal torus $T$ of $B_1$ such that $B_1,T$ are normalized by $u$. Then 
$B_2=B_1U\in\cb$ is also normalized by $u$. Let $G^0_{der}$ be the derived group of $G^0$. Let $Y_J$ (resp. $Y$) 
be the group of $1$-parameter subgroups of $(\cz_L\cap G^0_{der})^0$ (resp. $T$). Let $\a_1,\do,\a_r$ be the 
simple roots $Y@>>>\ZZ$ of $G^0$ (relative to $B_2,T$) which are not simple roots of $L$ (relative to $B_1,T$). 
Now $\QQ\ot Y_J$ has a unique basis $b_1,\do,b_r$ such that $\a_i(b_j)=\d_{ij}$ for all $i,j$. Let $d\in\ZZ_{>0}$
be such that $db_1,\do,db_r$ belong to $Y_J$. Let $\l=db_1+\do+db_r\in Y_J$ (in additive notation). Then 
$t\m\Ad(\l(t))$ has $>0$ weights on the Lie algebra of $U$. It remains to show that $\l(t)$ commutes with $u$. 
This follows from the fact that $\Ad(u)$ in its natural action on $Y_J$ permutes among themselves the elements 
$db_1,\do,db_r$. This completes the proof of (b) hence that of 1.3(a).

Let $\g\in\uD$. We can find $J\sub S$ such that $\e_D(J)=J$ and $\g'\in\uD_J$ such that $\g'$ is distinguished 
relative to $N_GL_J,D_J$. Using 1.3(d) we can find $C'\in(\uWW_J)_{D_J}^{el}$ and $\g'=\Ph(C')$ where $\Ph$ is 
defined relative to $N_GL_J,D_J$. Let $C$ be the $\e_D$-conjugacy class of $\WW$ that contains $C'$. By an 
earlier part of the argument, $\Ph(C)$ is well defined (relative to $G,D$) and is the unique $G^0$-conjugacy 
class that contains $\g'$; since $\g$ has the same property we have $\Ph(C)=\g$. This proves 1.3(c). This 
completes the proof of (a).

\subhead 1.5\endsubhead
Clearly, to prove that 1.3 holds for $G,D$ we may replace $G$ by the subgroup generated by $D$ hence we may 
assume that

(a) {\it $G/G^0$ is cyclic with generator $D$.}
\nl
Until the end of 1.11 we assume that (a) holds. In the case where $D=G^0$ (so that $G=G^0$), 1.3 follows from 
\cite{\WE, 0.4, 0.6} when $p$ is not a bad prime for $G$ and from the Addendum at the end of \cite{\WE} (based on
\cite{\LX} and \cite{\WEH}) when $p$ is a bad prime for $G$. We now consider the general case. In 1.6-1.11 we 
will give a number of reductions of the theorem.

\subhead 1.6\endsubhead
Let $G'=G/\cz_{G^0}^0$. Let $D'=\p(D)$ (a connected component of $G'$). We may identify $\WW$ with the Weyl group
of $G'{}^0$, $\e_{D'}:\WW@>>>\WW$ with $\e_D$ and $\WW_{D-min}$ with $\WW_{D'-min}$. We may identify $\uD'$ 
(defined in terms of $G',D'$) with $\uD$ via $\g\m\g'=\{g\in\p\i(\g);g\text{ unipotent}\}$, see 
\cite{\CDII, 12.2(a)}; this also identifies $\uD'_{dis}$ with $\uD_{dis}$ (see \cite{\CDII, 12.2(b)} and its 
proof). Then for $w\in\WW_{D-min}=\WW_{D'-min}$ we have $\Si_{w,D}=\Si_{w,D'}$ as partially ordered sets. Hence 
1.3 holds for $G',D'$ if and only if it holds for $G,D$. 

\subhead 1.7\endsubhead
Assume that $G$ is such that $G^0$ is semisimple. We can find a reductive group $\ti G$ with $\ti G^0$ 
semisimple, simply connected and a surjective homomorphism of algebraic groups $\p:\ti G@>>>G$ such that
$\ker\p\sub\cz_{\ti G^0}$. Then $\ti G^0=\p\i(G^0)$ and $\ti D=\p\i(D)$ is a connected component of $\ti G$.
Moreover, the obvious map $\ti G/\ti G^0@>>>G/G^0$ is a bijection carrying $\ti D$ to $D$. We may identify $\WW$ 
with the Weyl group of $\ti G^0$, $\e_{\ti D}:\WW@>>>\WW$ with $\e_D$ and $\WW_{D-min}$ with $\WW_{\ti D-min}$. 
We may identify $\un{\ti D}$ (defined in terms of $\ti G,\ti D$) with $\uD$ via 
$\g\m\g'=\{g\in\p\i(\g);g\text{ unipotent}\}$, see \cite{\CDII, 12.3(a)}. This also identifies 
$\un{\ti D}_{dis}$ with $\uD_{dis}$ (see \cite{\CDII, 12.3(b)} and its proof). Then for 
$w\in\WW_{D-min}=\WW_{\ti D-min}$ we have $\Si_{w,D}=\Si_{w,\ti D}$ as partially ordered sets. Hence if 1.3 holds
for $\ti G,\ti D$ then it holds for $G,D$.

\subhead 1.8\endsubhead
Next we assume that $G$ is such that $G^0$ is semisimple, simply connected. We can write uniquely $G^0$ as a 
product $G^0=G_1\T G_2\T\do\T G_k$ where each $G_i$ is a closed connected normal subgroup of $G$ different from 
$\{1\}$ and minimal with these properties. For $i\in[1,k]$ let 
$G'_i=G/(G_1\T\do\T G_{i-1}\T G_{i+1}\T\do\T G_k)$. Then $G'_i$ is an affine algebraic group with $G'_i{}^0=G_i$ 
and the image of $D$ in $G'_i$ is a connected component $D_i$ of $G'_i$. Also we have an obvious imbedding of
algebraic groups $G@>>>G'_1\T G'_2\T\do\T G'_k$ by which we identify $G$ with a closed subgroup of 
$G'_1\T G'_2\T\do\T G'_k$ with the same identity component; then $D$ becomes $D_1\T D_2\T\do\T D_k$. From the 
definitions we have a natural bijection
$$\uD_1\T\uD_2\T\do\T\uD_k@>>>\uD,\qua (\g_1,\g_2,\do,\g_k)\m\g_1\T\g_2\T\do\T\g_k$$
which restricts to a bijection 
$$(\uD_1)_{dis}\T(\uD_2)_{dis}\T\do\T(\uD_k)_{dis}@>>>\uD_{dis}.$$
Let $\WW_i$ be the Weyl group of $G_i$. Then $(\WW_i)_{D_i-min}$ is defined in terms of $G'_i,D_i$ in the same
way as $\WW_{D-min}$ is defined in terms of $G,D$ and we have canonically (compare \cite{\GP, Exercise 3.10} in 
the case where $D=G^0$): 
$$(\WW_1)_{D_1-min}\T(\WW_2)_{D_2-min}\T\do\T(\WW_k)_{D_k-min}@>\si>>\WW_{D-min}.$$
If $(w_1,w_2,\do,w_k)\lra w$ under the last bijection then from the definition we can identify
$\Si_{w_1,D_1}\T\Si_{w_2,D_2}\T\do\T\Si_{w_k,D_k}$ with $\Si_{w,D}$ as partially ordered sets. Hence if 1.3 holds
for each of $G'_i,D_i$ then it holds for $G,D$.

\subhead 1.9\endsubhead
Next we assume that $G$ is such that $G^0$ is semisimple, simply connected, $\ne\{1\}$ and that $G$ has no closed
connected normal subgroups other than $G^0$ and $\{1\}$. We have $G^0=H_0\T H_1\T\do\T H_{m-1}$ where $H_i$ are
connected, simply connected, almost simple, closed subgroups of $G^0$. Let $c\in D$ be a unipotent
quasi-semisimple element (see \cite{\CDI, 1.4, 1.9}). We can assume that $H_i=c^iH_0c^{-i}$ for $i\in[0,m-1]$ and
$c^mH_0c^{-m}=H_0$. Let $G'$ be the subgroup of $G$ generated by $H_0$ and $c^m$. Then $G'$ is closed,
$G'{}^0=H_0$ and $D'=c^mH_0$ is a connected component of $G'$. 

Define $D'@<a<<G^0\T D'@>b>>D$ by $a(g,c^mh)=c^mh$, $b(g,c^mh)=gchg\i$ (with $h\in H_0$); we have a bijection 
$\t:\uD'@>\si>>\uD$ given by $\g'\m\g$ where $a\i(\g')=b\i(\g)$, see \cite{\CDII, 12.5(a),(c)}. This restricts to
a bijection $\uD'_{dis}@>\si>>\uD_{dis}$ (see \cite{\CDII, 12.5(b)}).

Let $\cb_i$ be the variety of Borel subgroups of $H_i$ ($i\in[0,m-1]$). Any $B\in\cb$ can be written uniquely in 
the form $B=B_0B_1\do B_{m-1}$ where $B_i\in\cb_i$ ($i\in[0,m-1]$). Let $\WW_i$ be the Weyl group of $H_i$ and 
let $\ul_i:\WW_i@>>>\NN$ be its length function. We can identify $\WW=\WW_0\T\WW_1\T\do\T\WW_{m-1}$ in an obvious
way so that $\ul(w_0,\do,w_{m-1})=\ul_0(w_0)+\do+\ul_{m-1}(w_{m-1})$ for $w_i\in\WW_i$. We have 
$\e_D(\WW_i)=\WW_{i+1}$ for $i\in[0,m-2]$, $\e_D(\WW_{m-1})=\WW_0$. Let $\e_{D'}:\WW_0@>>>\WW_0$ be the 
automorphism defined by $D'$. We have $\e_{D'}(v)=\e_D^m(v)$ for $v\in\WW_0$.

In this subsection we assume that 1.3 holds for $G',D'$ and we show that it then also holds for $G,D$.

Let $C\in\uWW_D$. To verify that 1.3(a) holds for $C$ we choose $v\in C_{min}$ such that $v\in(\WW_0)_{D'-min}$ 
(see \cite{\GKP, 2.7}). Let $C'$ be the $\e_{D'}$-conjugacy class of $v$ in $\WW_0$. Let $\g'\in\uD',\g\in\uD$ be
such that $\g=\t(\g')$. We show that:

(a) {\it $\g'\in\Si_{v,D'}$ if and only if $\g\in\Si_{v,D}$.}
\nl
Let $h\in H_0$ be such that $c^mh\in\g'$; then $ch\in\g$. It is enough to show that the sets 
$$Z=\{B_0\in\cb_0;(B_0,c^mhB_0(c^mh)\i)\in\co_{v;H_0}\},$$
$$\align Z'=&\{(B_0,B_1,\do,B_{m-1})\in\cb_0\T\cb_1\t\cb_{m-1};(B_0,chB_{m-1}(ch)\i)\in\co_{v,H_0},\\&
B_1=chB_0(ch)\i,\do,B_{m-1}=chB_{m-2}(ch)\i\}\endalign$$
are in bijection (here $\co_{v,H_0}$ is defined like $\co_v$ but relative to $H_0$ instead of $G^0$). Now $Z'$ is
clearly in bijection with $\{B_0\in\cb_0;(B_0,(ch)^mB_0(ch)^{-m})\in\co_{v,H_0}\}$. It is enough to show that for 
$B_0\in\cb_0$ we have $(ch)^mB_0(ch)^{-m}=c^mhB_0(c^mh)\i$. We have $(ch)^m=c^mzh=c^mhz$ where 
$$z=(c^{-m+1}hc^{m-1})(c^{-m+2}hc^{m-2})\do(c^{-1}hc)\in H_1H_2\do H_{m-1}$$ 
commutes with $H_0$ hence $zB_0z\i=B_0$ and hence 
$$(ch)^mB_0(ch)^{-m}=c^mhzB_0z\i(c^mh)\i=c^mhB_0(c^mh)\i,$$
as required. This proves (a).

From (a) it follows that $\t$ defines a bijection $\Si_{v,D'}@>\si>>\Si_{v,D}$. Now $\t$ is compatible with the 
partial order (this follows from the proof of \cite{\CDII, 12.5(a)}). Since 1.3(a) holds for $G',D',C'$, it also
holds for $G,D,C$. Thus 1.3(a) holds for $G,D$. Similarly, 1.3(b)-(d) hold for $G,D$. (Note that $v$ above is 
$\e_D$-elliptic in $\WW$ if and only if it is $\e_{D'}$-elliptic in $\WW_0$.)

\subhead 1.10\endsubhead
Next we assume that $G$ is such that $G^0$ is semisimple, simply connected, almost simple. Let $\D$ be the 
subgroup of $\cz_G$ consisting of all unipotent elements in $\cz_G$. Let $G'=G/\D$ and let $\p:G@>>>G'$ be the 
obvious homomorphism. Now $\p$ induces an isomorphism $G^0@>\si>>G'{}^0$ and we have $\cz_{G'}\sub G'{}^0$, see 
\cite{\CDII, 12.6(a)}). Let $D'=\p(D)$, a connected component of $G'$. Then $\p$ restricts to an isomorphism 
$D@>\si>>D'$. We may identify $\WW$ with the Weyl group of $G'{}^0$, $\e_{D'}:\WW@>>>\WW$ with $\e_D$ and 
$\WW_{D-min}$ with $\WW_{D'-min}$. We may identify $\uD'$ (defined in terms of $G',D'$) with $\uD$ via 
$\g'\m\g=D\cap\p\i(\g')$; this also identifies $\uD'_{dis}$ with $\uD_{dis}$ (see \cite{\CDII, 12.6}). Then for 
$w\in\WW_{D-min}=\WW_{D'-min}$ we have $\Si_{w,D}=\Si_{w,D'}$ as partially ordered sets. Hence if 1.3 holds for 
$G',D'$, then it holds for $G,D$.

\subhead 1.11\endsubhead
We now assume that $G$ is such that $G^0$ is semisimple, almost simple and $\cz_G\sub G^0$. Let $G'=G/\cz_{G^0}$ 
and let $\p:G@>>>G'$ be the obvious homomorphism. Let $D'=\p(D)$, a connected component of $G'$. We may identify 
$\WW$ with the Weyl group of $G'{}^0$, $\e_{D'}:\WW@>>>\WW$ with $\e_D$ and $\WW_{D-min}$ with $\WW_{D'-min}$. We
may identify $\uD'$ (defined in terms of $G',D'$) with $\uD$ via $\g'\m\g=\p\i(\g')$; when $G=G^0$ this is 
obvious, when $G\ne G^0$, see \cite{\CDII, 12.7(b)}. This also identifies $\uD'_{dis}$ with $\uD_{dis}$ (see
\cite{\CDII, 12.7(c)}). Then for $w\in\WW_{D-min}=\WW_{D'-min}$ we have $\Si_{w,D}=\Si_{w,D'}$ as partially 
ordered sets. Hence 1.3 holds for $G',D'$ if and only if it holds for $G,D$.

\subhead 1.12\endsubhead
We now discuss the proof of Theorem 1.3. If $G^0=\{1\}$, the result is trivial. We can assume that $\dim(G^0)>0$ 
and that 1.3 is already known for any $G',D'$ with $\dim(G'{}^0)<\dim(G^0)$. From the results in 1.6-1.11 we see 
that we may assume that $G^0$ is semisimple, adjoint, almost simple, with $\cz_G\sub G^0$ and $D$ generates $G$. 
Moreover, as we have seen in 1.5, we can assume that $D\ne G^0$. Then, as in \cite{\CDII, 12.7}, we must be in 
one of the following four cases:

(a) $G^0=PGL_m(\kk)$, $m\ge3$, $p=2$;

(b) $G^0=PSO_{2m}(\kk)$, $m\ge4$, $p=2$;

(c) $G^0=PSO_8(\kk)$, $p=3$;

(d) $G^0$ is adjoint of type $E_6$, $p=2$;
\nl
moreover, $|G/G^0|=p$ and conjugation by an element of $D$ is not an inner automorphism of $G^0$. For $G,D$ as in
(c) or (d), the proof of 1.3 is given in \S2. If $G,D$ are as in (b), we see from 1.11 that it is enough to prove
1.3 when $G,D$ are replaced by 

(b${}'$) $G'=O_{2m}(\kk),D'=G'-G'{}^0$, $m\ge4,p=2$; 
\nl
the proof of 1.4 in this case is given in \S3; then 1.3 holds in this case by 1.4(a). If $G,D$ are as in (a), we 
see from 1.6 that it is enough to prove 1.3 when $G,D$ are replaced by $G',D'$ with 

(a${}'$) $G'{}^0=GL_m(\kk)$, $|G'/G'{}^0|=2$, some element of $D'$ acts on $G'{}^0$ by conjugation as 
$(a_{ij})\m(a_{ji})\i$, $m\ge3,p=2$; 
\nl
the proof of 1.4 in this case is given in \S5 (based on results in \S4); then 1.3 holds in this case by 1.4(a).

\subhead 1.13\endsubhead
In the remainder of this section, $G,D$ are as in 1.1(a),(b).

\proclaim{Theorem 1.14} Let $C\in\uWW_D^{el}$, $w\in C_{min}$ and let $\g=\Ph(C)$, see 1.3(a). Let $g\in\g$. We 
have $\dim(Z_{G^0}(g)/\cz_{G^0}^D)=\ul(w)$.
\endproclaim
We first go through a sequence of reductions as in 1.6-1.11.

Assume first that we are in the setup of 1.6 and that the theorem holds for $G',D'$. We have $\cz_{G'{}^0}=\{1\}$
and it is enough to show that $Z_{G^0}(g)/\cz_{G^0}^D\cong Z_{G'{}^0}(\p(g))$. This follows from 
\cite{\CDII, 12.2(b)}.

Assume now that we are in the setup of 1.7 and that the theorem holds for $\ti G,\ti D$. Let $\tg$ be a unipotent
element in $\p\i(g)$. It is enough to show that $\dim (Z_{\ti G^0}(\tg))=\dim(Z_{G^0}(g)$. This follows from 
\cite{\CDII, 12.3(b)}.

Assume now that we are in the setup of 1.8 and that the theorem holds for $G'_i,D_i$ ($i\in[1,k]$). Then clearly 
the theorem holds for $G,D$. 

Assume now that we are in the setup of 1.9 and that the theorem holds for $G',D'$. It is enough to show that if 
$h\in H_0$ then $Z_{G^0}(ch)=Z_{H_0}(c^mh)$. This follows from \cite{\CDII, 12.5(b)}.

Assume now that we are in the setup of 1.10 and that the theorem holds for $G',D'$. It is enough to show that
$Z_{G^0}(g)\cong Z_{G'{}^0}(\p(g))$. This is shown in \cite{\CDII, 12.6}.

By the arguments above the proof of the theorem is reduced to the special case where $G$ is as in 1.11; we can 
also assume that $D$ generates $G$. If $D=G^0$ then the theorem is already known, see \cite{\WE, 4.4(b)}. Thus we
can assume in addition that $D\ne G^0$. Assume now that (in the setup of 1.11), the theorem holds for $G',D'$. 
Using \cite{\CDII, 12.7(c)} we deduce that the theorem holds for $G,D$. We see that it is enough to prove the
theorem in the cases 1.12(a)-(d). If we are in the case 1.12(c) or 1.12(d) the result follows from the explicit 
description of the map $\Ph$ in \S2 (the values of $\ul(w)$ can be extracted from the character table of the 
appropriate Hecke algebra available through the CHEVIE package). The proof in the case 1.12(b) (or equivalently 
1.12(b${}'$)) is almost identical to the proof for $G^0$ given in \cite{\WE, 4.4(b)} and will be omitted. The 
case 1.12(a) (or equivalently 1.12(a${}'$)) is treated in 5.10. 

\proclaim{Theorem 1.15} Let $C\in\uWW_D^{el}$, $w\in C_{min}$ and let $\g=\Ph(C)$, see 1.3(a). The $G^0$-action 
$x:(g,B)\m(xgx\i,xBx\i)$ on $\fB_w^\g:=\{(g,B)\in\g\T\cb;(B,gBg\i)\in\co_w\}$ is transitive.
\endproclaim
We set $\fB_w^D=\{(g,B)\in D\T\cb;(B,gBg\i)\in\co_w\}$. We state the following result which is similar to 
\cite{\WE, 5.2(a)}. 

(a) {\it If $w',w''\in C_{min}$ then there exists an isomorphism $\fB_{w'}^D@>\si>>\fB_{w''}^D$ commuting with 
the $G^0$-actions and commuting with the first projections $\fB_{w'}^D@>>>D$, $\fB_{w''}^D@>>>D$.}
\nl
The proof of (a) is along the same lines as that in \cite{\WE, 5.3}. Using a result of \cite{\GP, 3.2.7} and its 
extension to the twisted case \cite{\GKP},\cite{\HE}, we see that we can assume that there exist
$b,c,b'\in\WW$ such that $w'=bc$, $w''=cb'$, $\e_D(b)=b'$, $\ul(b)+\ul(c)=\ul(bc)=\ul(cb')$. If 
$(g,B)\in\fB_{bc}^D$ then there is a unique $B'\in\cb$ such that $(B,B')\in\co_b$, $(B',gBg\i)\in\co_c$. We have
$(gBg\i,gB'g\i)\in\co_{b'}$ hence $(B',gB'g\i)\in\co_{cb'}$ so that $(g,B')\in\fB_{cb'}^D$. Thus we have defined 
a morphism $\a:\fB_{bc}^D@>>>\fB_{cb'}^D$, $(g,B)\m(g,B')$. Similarly, if $(g,B')\in\fB_{cb'}^D$, there exists a 
unique $B''\in\cb$ such that $(B',B'')\in\co_c$, $(B'',gB'g\i)\in\co_{b'}$; we have $(g,g\i B''g)\in\fB_{bc}^D$. 
Thus we have defined a morphism $\a':\fB_{cb'}^D@>>>\fB_{bc}^D$, $(g,B')\m(g,g\i B''g)$. From the definition it 
is clear that $\a'\a(g,B)=(g,B)$ for all $(g,B)\in\fB_{bc}^D$ and $\a\a'(g,B')=(g,B')$ for all 
$(g,B')\in\fB_{cb'}^D$. It follows that $\a,\a'$ are isomorphisms. They have the required properties.

From (a) we see that:

(b) {\it If the theorem is true for some $w\in C_{min}$ then it is true for any $w\in C_{min}$.}
\nl
Note that the following is an equivalent formulation of the theorem.

(c) {\it In the setup of the theorem let $B\in\cb$ and let $g,g'\in\g$ be such that $(B,gBg\i)\in\co_w$, 
$(B,g'Bg'{}\i)\in\co_w$. Then there exists $x\in B$ such that $xgx\i=g'$.} 
\nl
Before proving the theorem (or equivalently (c)) we go through a sequence of reductions as in 1.6-1.11.

Assume first that we are in the setup of 1.6 and that (c) holds for $G',D'$. If $B,g,g'$ are as in (c), we have 
$\p(g)\in\p(\g),\p(g')\in\p(\g)$ and $\p(B)$ is a Borel subgroup of $G'$. We can find $x'\in\p(B)$ such that 
$x'\p(g)x'{}\i=\p(g')$. We have $x'=\p(x)$ where $x\in B$ and $xgx\i=zg'$ for some $z\in\cz_{G^0}^0$. Using 
\cite{\CDI, 1.3(a)} we can write $z=ytg'y\i g'{}\i$ with $t,y\in\cz_{G^0}^0$, $tg'=g't$. We have 
$zg'=ytg'y\i=xgx\i$, $(x\i y)tg'(x\i y)\i=g$. In particular $tg'$ is unipotent. Since $g'$ is unipotent and $t$ 
is semisimple and commutes with $g'$ it follows that $t=1$. Hence $(x\i y)g'(x\i y)\i=g$. We see that (c) holds 
for $G,D$. On the other hand if (c) holds for $G,D$ then it obviously holds for $G',D'$.

Assume now that we are in the setup of 1.7 and that (c) holds for $\ti G,\ti D$. Then clearly (c) holds for
$G,D$. 

Assume now that we are in the setup of 1.8 and that (c) holds for $G'_i,D_i$ ($i\in[1,k]$). Then clearly (c) 
holds for $G,D$. 

Assume now that we are in the setup of 1.9 and that the theorem holds for $G',D'$. We show that the theorem holds
for $G,D$. By (b) we can assume that $w=v\in C_{min}\cap(\WW_0)_{D'-min}$. Let $(g,B)\in\fB_v^\g$, 
$(g',B')\in\fB_v^\g$. It is enough to show that $(g,B),(g',B')$ are in the same $G^0$-orbit. Replacing 
$(g,B),(g',B')$ by pairs in the same $G^0$-orbit we can assume that $g=ch_0,g'=ch'_0$ where $h_0,h'_0\in H_0$ and
$c^mh_0,c^mh'_0$ are in the same (unipotent) $H^0$-conjugacy class $\g'$ in $D'$. We write $B=B_0B_1\do B_{m-1}$,
$B'=B'_0B'_1\do B'_{m-1}$ with $B_i,B'_i\in\cb_i$. We have 
$$(B_0,ch_0B_{m-1}(ch_0)\i)\in\co_{v,H_0}, B_1=ch_0B_0(ch_0)\i,\do,B_{m-1}=ch_0B_{m-2}(ch_0)\i,$$
$$(B'_0,ch'_0B'_{m-1}(ch'_0)\i)\in\co_{v,H_0}, B'_1=ch'_0B'_0(ch'_0)\i,\do,B'_{m-1}=ch'_0B'_{m-2}(ch'_0)\i.$$
As in 1.9 from this we deduce that $(B_0,c^mh_0B_0(c^mh_0)\i)\in\co_{v;H_0}$ and similarly
$(B'_0,c^mh'_0B_0(c^mh'_0)\i)\in\co_{v;H_0}$. Since the theorem holds for $G',D'$ we can find $x_0\in H_0$ such 
that $B'_0=x_0B_0x_0\i$, $c^mh'_0=x_0c^mh_0x_0\i$. For $i\in[1,m-1]$ we set $x_i=c^ih'_0x_0h_0\i c^{-i}\in H_i$; 
we show that $x_iB_ix_i\i=B'_i$. An equivalent statement is: 
$$c^ih'_0x_0h_0\i c^{-i}(ch_0)^iB_0(ch_0)^{-i}c^ih_0x_0\i h'_0{}\i c^{-i}=(ch'_0)^ix_0B_0x_0\i(ch'_0)^{-i}$$
It is enough to show that if $\l$ is defined by 
$$c^ih'_0x_0h_0\i c^{-i}(ch_0)^i=(ch'_0)^ix_0\l$$
then $\l\in H_1H_2\do H_{m-1}$ so that $\l B_0\l\i=B_0$. We have 
$$\l=x_0\i(ch'_0)^{-i}c^ih'_0x_0h_0\i c^{-i}(ch_0)^i=x_0\i z'x_0z$$
where $z=h_0\i c^{-i}(ch_0)^i$ and $z'\in (ch'_0)^{-i}c^ih'_0$ belong to $H_1H_2\do H_{m-1}$ hence they
commute with $H_0$. Thus $\l=z'z\in H_1H_2\do H_{m-1}$.

Now let $x=x_0x_1\do x_{m-1}\in G^0$. We have 
$$xBx\i=(x_0B_0x_0\i)(x_1B_1x_1\i)\do(x_{m-1}B_{m-1}x_{m-1}\i)=B'_0B'_1\do B'_{m-1}=B'.$$
We show that $xgx\i=g'$ that is, $x_0x_1\do x_{m-1}=ch'_0x_0x_1\do x_{m-1}h_0\i c\i$. An equivalent statement is 
that
$$x_0x_1\do x_{m-1}=(cx_{m-1}c\i)(ch'_0x_0h_0c\i)(cx_1c\i)\do(cx_{m-2}c\i).$$
(We use the fact that $H_0,H_1,\do,H_{m-1}$ commute with each other.) It is enough to show that
$$x_0=cx_{m-1}c\i,x_1=ch'_0x_0h_0c\i, x_2=cx_1c\i,\do,x_{m-1}=cx_{m-2}c\i.$$
These equalities, except the first, follow from the definition of $x_1,\do,x_{m-1}$; the first equality is the 
same as $x_0=c^mh'_0x_0h_0\i c^{-m}$; it holds by the definition of $x_0$. Thus, the theorem holds for $G',D'$. 

Assume now that we are in the setup of 1.10 and that (c) holds for $G',D'$. If $B,g,g'$ are as in (c), we have 
$\p(g)\in\p(\g),\p(g')\in\p(\g)$ and $\p(B)$ is a Borel subgroup of $G'$. We can find $x'\in\p(B)$ such that 
$x'\p(g)x'{}\i=\p(g')$. We have $x'=\p(x)$ where $x\in B$ and $xgx\i=zg'$ for some $z\in\D$ such that the 
connected component containing of $zg'$ (that is $zD$) is equal to the connected component containing $xgx\i$ 
(that is $D$). We see that $zD=D$ hence $z\in G^0\cap\G=\{1\}$. Thus $xgx\i=g'$ and (c) holds for $G,D$.

By the arguments above the proof of the theorem is reduced to the special case where $G$ is as in 1.11; we can 
also assume that $D$ generates $G$. If $D=G^0$ then the theorem is already known, see \cite{\WEH, 0.2}. Thus we 
can assume in addition that $D\ne G^0$. Then, as in \cite{\CDII, 12.7(d)}, we have that for any $g'\in D$, the 
homomorphism $\cz_{G^0}@>>>\cz_{G_0}$, $y\m g'{}\i yg'y\i$ is an isomorphism. 

Assume now that (in the setup of 1.11), (c) holds for $G',D'$. If $B,g,g'$ are as in (c), we have 
$\p(g)\in\p(\g),\p(g')\in\p(\g)$ and $\p(B)$ is a Borel subgroup of $G'$. We can find $x'\in\p(B)$ such that 
$x'\p(g)x'{}\i=\p(g')$. We have $x'=\p(x)$ where $x\in B$ and $xgx\i=g'z$ for some $z\in\cz_{G^0}$. As noted 
above, we can write $z=g'{}\i yg'y\i$ with $y\in\cz_{G^0}$. Then $g'z=yg'y\i$ and $xgx\i=yg'y\i$ so that 
$(y\i x)g(y\i x)\i=g'$; note also that $y\i x\in B$. Thus (c) holds for $G,D$. 

We see that it is enough to prove (c) in the cases 1.12(a)-(d). In the cases 1.12(c),(d) the theorem is contained
in 2.3(b), 2.4(b). The proof in the case 1.12(b) (or equivalently 1.12(b${}'$)) is almost identical to the proof 
for $G^0$ given in \cite{\WEH, \S3} and will be omitted. The case 1.12(a) (or equivalently 1.12(a${}'$)) is 
treated in 5.12. 

\proclaim{Theorem 1.16} For any $\g\in\uD$ the function $\Ph\i(\g)@>>>\NN$, $C\m \mu(C)$ ($\mu$ as in the proof 
of 1.4(a)) reaches its minimum at a unique element $C_0\in\Ph\i(\g)$. Thus we have a well defined map 
$\Ps:\uD@>>>\uWW_D$, $\g\m C_0$ such that $\Ph\Ps:\uD@>>>\uD$ is the identity map.
\endproclaim
By a sequence of reductions as in 1.6-1.11 we see that it is enough to prove the theorem assuming that $D=G^0$ or
that we are in one of the cases 1.12(a)-(d). If $D=G^0$ the theorem follows from \cite{\WEII, 0.2}. If we are in
the case 1.12(c),(d), the result follows from the tables in 2.3, 2.4. If we are in the case 1.12(a) or (b), or
equivalently in the case 1.12(a${}'$) or 1.12(b${}'$), the result follows by arguments similar to those in 
\cite{\WEII}.

\head 2. The exceptional cases\endhead
\subhead 2.1\endsubhead
In this section we assume that $G,D$ are as in 1.12(a)-(d) and that 
$$\kk,F_q,\cf,T_w:\cf@>>>\cf,\ch_q,S_x:\cf@>>>\cf$$
are as in the proof of 1.2(a). Let $u$ be a unipotent quasi-semisimple element of $D(F_q)$. Then
$u$ has order $p$ and $\e_D^p=1:\WW@>>>\WW$. Let $E^i(i\in I)$ be a set of representatives for the isomorphism 
classes of irreducible representations of $\WW$ over $\bbq$. For $i\in I$ let $E^i_q$ be the irreducible 
$\ch_q$-module corresponding to $E^i$ (it depends on a fixed choice of $\sqrt{q}$ in $\bbq$) and let 
$\r_i=\Hom_{\ch_q}(E^i_q,\cf)$. We regard $\r_i$ as an (irreducible) $G^0(F_q)$-module by $x:\ph\m\ph'$ 
($x\in G^0(F_q)$) where $\ph'(e)=S_x\ph(e)$ for $e\in E^i_q$. We have an isomorphism 
$$A:\op_{i\in I}\r_i\ot E^i_q@>\si>>\cf\tag a$$
given by $\ph\ot e\m\ph(e)$ for $\ph\in\r_i,e\in E^i_q$. Let $I_{ex}$ be the set of all $i\in I$ such that there 
exists a linear map $v_i:E^i@>>>E^i$ with $v_iw=\e_D(w)v_i:E^i@>>>E^i$ for all $w\in\WW$, $v_i^p=1$. Note that 
for $i\in I_{ex}$, $v_i$ is only defined up to multiplication by a $p$-th root of $1$. However there is a 
canonical choice for $v_i$, the "preferred extension", see \cite{\CSIV, 17.2}; we shall assume that $v_i$ is this
canonical choice. As in \cite{\OR, p.61}, $v_i$ gives rise to a linear map $V_i:E^i_q@>>>E^i_q$ such that 
$V_iT_w=T_{\e_D(w)}V_i:E^i_q@>>>E^i_q$ for all $w\in\WW$ and $V_i^p=1$. 

For $i\in I_{ex}$ we define a linear map $U_i:\r_i@>>>\r_i$ by $\ph\m\ph'$ where $\ph'(e)=S_u\ph(V_i\i(e))$ for 
$e\in E^i_q$. (This is well defined since $T_{\e_D(y)}S_g=S_gT_y:\cf@>>>\cf$ for any $y\in\WW,g\in D$, see 1.2.) 
We have $U_i^c=1$ and $U_ix=(uxu\i)U_i:\r_i@>>>\r_i$ for any $x\in G^0(F_q)$. Hence we can regard $\r_i$ as a 
$G(F_q)$-module extending the $G^0(F_q)$-module considered above so that $u$ acts as $U_i$; we call this the
{\it preferred extension} of the $G^0(F_q)$-module $\r_i$ to a $G(F_q)$-module.

Let $g\in D(F_q),w\in\WW$. We write $g=ux,x\in G^0(F_q)$. For any $i\in I_{ex}$ we define a linear map 
$H_i:\r_i\ot E^i_q@>>>\r_i\ot E^i_q$ by $H_i(\ph\ot e)=(g\ph)\ot(V_iT_w(e))$. For $\ph\in\r_i,e\in E^i_q$ we have 
$S_gT_w(A(\ph\ot e))=A(H_i(\ph\ot e))$. (The left hand side is $S_xS_uT_w(\ph(e))=S_xS_u(\ph(T_we))$; the right 
hand side is $(xU_i\ph)(V_iT_we)=S_xS_u\ph(V_i\i V_iT_we)=S_xS_u(\ph(T_we))$, as desired. Thus the endomorphism 
$S_gT_w$ of $\cf$ corresponds under the isomorphism (a) to an endomorphism $R$ of $\op_{i\in I}\r_i\ot E^i_q$ 
such that $R|_{\r_i\ot E^i_q}=H_i$ (if $i\in I_{ex}$; it is clear that $R$ permutes the summands $\r_i\ot E^i_q$ 
with $i\in I-I_{ex}$ according to a fixed point free permutation of $I-I_{ex}$. It follows that
$$\tr(S_gT_w:\cf@>>>\cf)=\sum_{i\in I_{ex}}\tr(H_i:\r_i\ot E^i_q@>>>\r_i\ot E^i_q)$$
hence
$$\tr(S_gT_w:\cf@>>>\cf)=\sum_{i\in I_{ex}}\tr(g,\r_i)\tr(V_iT_w,E^i_q).\tag b$$
Applying (b) with $g=u$, $w=1$ we obtain
$$\tr(S_u:\cf@>>>\cf)=\sum_{i\in I_{ex}}\tr(u,\r_i)\tr(V_i,E^i_q).\tag c$$

\subhead 2.2 \endsubhead
Let $F:G@>>>G$ be the Frobenius map corresponding to the $F_q$-rational structure of $G$. The induced map
$\cb@>>>\cb$ is denoted again by $F$. Let $\WW_0$ be the fixed point set of $\e_D:\WW@>>>\WW$. This is a subgroup
of $\WW$, in fact a Weyl group with generators indexed by the orbits of $\e_D:S@>>>S$.

For any $w\in\WW$ let $X_w=\{B\in\cb;(B,F(B))\in\co_w\}$, see \cite{\DL}. If $g\in D(F_q)$ and $B\in X_w$ then 
$(gBg\i,F(gBg\i))=(gBg\i,gF(B)g\i)\in\co_{\e_D(w)}$. Hence if $w\in\WW_0$, then $G(F_q)$ acts by conjugation on 
$X_w$; hence $G(F_q)$ acts on the $l$-adic cohomology with compact support $H^i_c(X_w,\bbq)$; we denote by $R_w$ 
the virtual representation $\sum_i(-1)^iH^i_c(X_w,\bbq)$ of $G(F_q)$. Let $\cg$ be the vector space of functions 
$D(F_q)@>>>\bbq$ generated by the characters of irreducible representations of $G^0(F_q)$ which appear in 
$H^i_c(X_w,\bbq)$ for some $w\in\WW_0,i\in\ZZ$ and are extendable to $G(F_q)$-modules. The character of $R_w$ 
($w\in\WW_0$) restricted to $D(F_q)$ belongs to $\cg$ and is denoted again by $R_w$. Also, if $i\in I_{ex}$, the 
character of the preferred extension $\r_i$ restricted to $D(F_q)$ belongs to $\cg$ and is denoted again by 
$\r_i$. For $r,r'$ in $\cg$ we set $(r,r')=|D(F_q)|\i\sum_{g\in D(F_q)}r(g)r'(g\i)$. Let $\cg_0$ be the vector 
space of all $f\in\cg$ such that $(f,\ti R_w)=0$ for all $w\in\WW_0$.

Let $\fe_j (j\in J)$ be a set of representatives for the irreducible representations of $\WW_0$. For $j\in J$ we 
define, following \cite{\LU} and \cite{\MA}, the element 
$$R_{\fe_j}=|\WW_0|\i\sum_{w\in\WW_0}\tr(w,\fe_j)R_w\in\cg.$$

\subhead 2.3\endsubhead
In this subsection we assume that we are in the setup of 1.12(d). In this case $I=I_{ex}$. Following \cite{\MA} 
we index the irreducible representations of $\WW$ as
$$\align&1_0,1_{36},10_9,6_1,6_{25},20_{10},15_5,15_{17},15_4,15_{16},20_2,20_{20},24_6,\\&
24_{12},30_3,30_{15},60_8,80_7,90_8,60_5,60_{11},64_4,64_{13},81_6,81_{10}.\endalign$$
This list is taken as the set $I$ so that $E^{1_0},\do,E^{81_{10}}$ in $\Irr(\WW)$, $E^{1_0}_q,\do,E^{81_{10}}_q$
(representations of $\ch_q$) and $\r_{1_0},\do,\r_{81_{10}}$ (representations of $G(F_q)$) are defined.

Following \cite{\MA}, we index the irreducible representations of $\WW_0$ (of type $F_4$) as 
$$\align&1_0,4_1,9_2,8'_3,8''_3,2'_4,2''_4,12_4,16_5,9'_6,9''_6,6'_6,6''_6,4'_7,4''_7,\\&
4_8,8'_9,8''_9,9_{10},1'_{12},1''_{12},4_{13},2'_{16},2''_{16},1_{24}.\endalign$$
This list is taken as the set $J$ so that $\fe_{1_0},\do,fe_{1_{24}}$ are in $\Irr(\WW_0)$; for $j\in J$ we write 
$R_j$ instead of $R_{\fe_j}$. From \cite{\MA, Theorem 8},\cite{\MAB} we see that for some function
$\e:I@>>>\{1,-1\}$ the following equalities (referred to as $\sp$) hold in $\cg$.
$$\e(1_0)\r_{1_0}=R_{1_0},$$
$$\e(1_{36})\r_{1_{36}}=R_{1_{24}},$$
$$\e(10_9)\r_{10_9}=(1/3)R_{12_4}+(2/3)R_{6'_6}-(1/3)R_{6''_6}+\Xi_1,$$
$$\e(6_1)\r_{6_1}=R_{2'_4},$$
$$\e(6_{25})\r_{6_{25}}=R_{2''_{16}},$$
$$\e(20_{10})\r_{20_{10}}=(1/6)R_{12_4}-(1/2)R_{4_8}+(1/3)R_{6'_6}+(1/3)R_{6''_6}-(1/2)R_{16_5}+\Xi_2,$$
$$\e(15_5)\r_{15_5}=(1/2)R_{2''_4}+(1/2)R_{9_2}-(1/2)R_{1'_{12}}-(1/2)R_{8'_3},$$
$$\e(15_{17})\r_{15_{17}}=(1/2)R_{2'_{16}}+(1/2)R_{9_{10}}-(1/2)R_{1''_{12}}-(1/2)R_{8''_9},$$
$$\e(15_4)\r_{15_4}=(1/2)R_{2''_4}-(1/2)R_{9_2}+(1/2)R_{1'_{12}}-(1/2)R_{8'_3},$$
$$\e(15_{16})\r_{15_{16}}=(1/2)R_{2'_{16}}-(1/2)R_{9_{10}}+(1/2)R_{1''_{12}}-(1/2)R_{8''_9},$$
$$\e(20_2)\r_{20_2}=R_{4_1},$$
$$\e(20_{20})\r_{20_{20}}=R_{4_{13}},$$
$$\e(24_6)\r_{24_6}=R_{8''_3},$$
$$\e(24_{12})\r_{24_{12}}=R_{8'_9},$$
$$\e(30_3)\r_{30_3}=(1/2)R_{2''_4}+(1/2)R_{9_2}+(1/2)R_{1'_{12}}+(1/2)R_{8'_3},$$
$$\e(30_{15})\r_{30_{15}}=(1/2)R_{2'_{16}}+(1/2)R_{9_{10}}+(1/2)R_{1''_{12}}+(1/2)R_{8''_9},$$
$$\e(60_8)\r_{60_8}=-(1/2)R_{12_4}+(1/2)R_{4_8}-(1/2)R_{16_5}+\Xi_3,$$
$$\e(80_7)\r_{80_7}=(1/6)R_{12_4}+(1/2)R_{4_8}+(1/3)R_{6'_6}+(1/3)R_{6''_6}+(1/2)R_{16_5}+\Xi_4,$$
$$\e(90_8)\r_{90_8}=(1/3)R_{12_4}-(1/3)R_{6'_6}+(2/3)R_{6''_6}+\Xi_5,$$
$$\e(60_5)\r_{60_5}=R_{4'_7},$$    
$$\e(60_{11})\r_{60_{11}}=R_{4''_7},$$  
$$\e(64_4)\r_{64_4}=\Xi_6,$$
$$\e(64_{13})\r_{64_{13}}=\Xi_7,$$
$$\e(81_6)\r_{81_6}=R_{9'_6},$$     
$$\e(81_{10})\r_{81_{10}}=R_{9''_6},$$   
where $\Xi_1,\Xi_2,\do,\Xi_7\in\cg_0$ and $\e(i)$ may apriori depend on $q$.
More precisely, in \cite{\MA}, for $i\in I$, the $G^0(F_q)$-module $\r_i$ is extended to a $G(F_q)$-module $\r'_i$
not by the preferred extension (2.1) but by the requirement that $(\r'_i,R_1)\ge0$ (this determines $\r'_i$ 
uniquely if $(\r'_i,R_1)>0$; if $(\r'_i,R_1)=0$ for one of the extensions then the same holds for the other 
extension and we pick arbitrarily one of the two extensions and call it $\r'_i$). Then the character of $\r'_i$ 
on $D(F_q)$ is an element of $\cg$ denoted again by $\r'_i$; we have $\r'_i=\pm\r_i$ for all $i\in I$. What is 
actually shown in \cite{\MA, Theorem 8},\cite{\MAB} is that the equations $\sp$ hold if $\e(i)\r_i$ is replaced 
by $\e'(i)\r'_i$ where $\e':I@>>>\{1,-1\}$ is given by $\e'(20_{10})=\e'(15_4)=\e'(15_{16})=\e'(60_8)=-1$, 
$\e'(i)=1$ for all other $i$. It follows that $\e(i)\r_i=\e'(i)\r'_i$ for all $i\in I-\{64_4,64_{13}\}$.

In \cite{\MAA}, the values $R_j(g)$ are explicitly computed (as polynomials in $q$) for any $j\in J$ and any 
unipotent element in $D(F_q)$. Since, by \cite{\MAA, Proposition 7}, we have $\sum_{g\in\g(F_q)}\Xi(g)=0$ for any
$\g\in\uD$, $\Xi\in\cg_0$, we see from $\sp$ that $\sum_{g\in\g(F_q)}\e(i)\r_i(g)$ is explicitly known (as a 
polynomial in $q$) for any $\g\in\uD$ and any $i\in I$.

In particular, if $\g\in\uD$ contains $u$ (see 2.1) then for each $i\in I$, $|\g(F_q)|\e(i)\r_i(u)$ is explicitly
known. Hence $\e(i)\r_i(u)$ is explicitly known as a polynomial in $q$. It turns out that the value of this
polynomial at $q=1$ is the integer $\tr(v_i,E^i)$. Using this and 2.1(c) we see that $\tr(S_u:\cf@>>>\cf)$ is a 
polynomial in $q$ whose value at $q=1$ is $\sum_{i\in I}\e(i)\tr(v_i,E^i)\tr(v_i,E^i)$. Now 
$\tr(S_u:\cf@>>>\cf)=|\{B\in\cb(F_q);uBu\i=B\}|$ is equal to the number of Borel subgroups defined over $F_q$ of 
a simple algebraic group of type $F4$ defined over $F_q$, hence is a polynomial in $q$ whose value at $q=1$ is 
$|WW_0|=1152$. Thus we have $\sum_{i\in I}\e(i)\tr(v_i,E^i)^2=1152$. By standard orthogonality relations for 
characters we have also $\sum_{i\in I}\tr(v_i,E^i)^2=1152$. Thus $\sum_{i\in I}(1-\e(i))\tr(v_i,E^i)^2=0$. Since 
$1-\e(i),\tr(v_i,E^i)^2$ are integers $\ge0$ it follows that for any $i\in I$ such that $\tr(v_i,E^i)\ne0$ we 
have $\e(i)=1$. Note that $\tr(v_i,E^i)=0$ for $i\in\{64_4,64_{13}\}$ and $\tr(v_i,E^i)\ne0$ for all other $i$. 
We see that $\e(i)=1$ for all $i\in I-\{64_4,64_{13}\}$. For $i\in\{64_4,64_{13}\}$ the equations in $\sp$ hold 
for any choice of $\e(i)$. Hence we may assume that $\e(i)=1$ for all $i\in I$. We now see that

(a) {\it $\sum_{g\in\g(F_q)}\r_i(g)$ is explicitly known (as a polynomial in $q$) for any $\g\in\uD$ and any 
$i\in I$.}
\nl
From 2.1(b) we have for any $\g\in\uD$ and any $w\in\WW_{D-min}$:
$$\align&|\fB_w^\g(F_q)|=|\{(g,B)\in\g(F_q)\T\cb(F_q);(B,gBg\i)\in\co_w\}|\\&=\sum_{g\in\g(F_q)}
\tr(S_gT_w:\cf@>>>\cf)=\sum_{i\in I}(\sum_{g\in\g(F_q)}\tr(g,\r_i))\tr(V_iT_w,E^i_q).\endalign$$
Here the last sum can be calculated as an explicit polynomial in $q$ using (a) and the known values of the
polynomials $\tr(V_iT_w,E^i_q)$ (available through the CHEVIE package). This calculation was performed using a 
computer. The explicit knowledge of the quantities $|\fB_w^\g(F_q)|$ as polynomials in $q$ allows us to decide 
when $\g\in\Si_{w,D}$ (this holds precisely when the polynomial is nonzero). We can then verify that Theorem 1.3 
holds in our case. (Note that the partial order on $\uD$ is known from \cite{\SP, p.250} and the description of 
$\uD_{dis}$ is known from \cite{\SP, p.161}.) The calculation of the sum above yields in particular:
$$|\fB_w^\g(F_q)|=|G^0(F_q)|$$
whenever $C\in\uWW^{el}$, $w\in C_{min}$ and $\g=\Ph(C)$. Now from \cite{\XW, 0.3(b)} we see that any $G^0$-orbit
on $\fB_w^\g$ (for the conjugation action on both factors) which is defined over $F_q$ has a number of 
$F_q$-rational points equal to $|G^0(F_q)|$; it follows that 

(b) {\it $\fB_w^\g$ is a single $G^0$-orbit if $C\in\uWW^{el}$, $w\in C_{min}$ and $\g=\Ph(C)$.}
\nl
We now describe explicitly the map $\Ph$ in our case. The $\e_D$-conjugacy classes in $\WW$ are in bijection with
the ordinary conjugacy classes in $\WW$ under the map $C\m Cw_0$ ($w_0$ is the longest element of $\WW$); we 
denote an $\e_D$-conjugacy class in $\WW$ by the same symbol (given in \cite{\CA}) as the corresponding ordinary 
conjugay class in $\WW$ (we also add an ${}^!$ to the name of an $\e_D$-elliptic conjugacy class). The objects of
$\uD$ are denoted as
$$\g_{52},\g_{36},\g_{30},\g_{28},\g_{24},\g_{22}^{20},\g_{22}^{14},\g_{18},\g_{16}^{14},\g_{16}^{12},\g_{14},
\g_{12},\g_{10}^8,\g_{10}^9,\g_8,\g_6,\g_4$$
where the subscript indicates the codimension of the class (it is taken from \cite{\SP} and 
\cite{\MAA, Table 10}) and the upperscript denotes the dimension of the largest unipotent subgroup of the
centralizer of an element in the class (it is taken from \cite{\MAA, Table 10}). The map $\Ph$ is as follows:

$2A_1,4A_1\m\g_{52}$
 
$A_0^!,A_1,3A_1\m\g_{36}$
 
$A_32A_1\m\g_{30}$

$A_3A_1\m\g_{24}$

$D_4\m \g_{28}$

$A_3\m\g_{22}^{20}$

$D_4(a_1)^!,D_5(a_1)\m\g_{18}$

$A_5A_1\m\g_{22}^{14}$

$E_6(a_2)^!,A_5\m\g_{16}^{14}$

$A_2^!,A_22A_1,A_2A_1\m\g_{16}^{12}$

$2A_2^!,2A_2A_1\m\g_{14}$

$3A_2^!\m\g_{12}\text{ dist}$   

$A_4A_1\m\g_{10}^8$

$D_5\m\g_{10}^9$

$A_4^!\m\g_8\text{ dist}$         

$E_6^!\m\g_6\text{ dist}$   

$E_6(a_1)^!\m\g_4\text{ dist}$.
\nl
Here we have indicated the distinguished unipotent classes by "dist".

\subhead 2.4\endsubhead
In this subsection we assume that we are in the setup of 1.12(c).  Following \cite{\MA} we index the irreducible
representations of $\WW$ in terms of 
$$I_{ex}=\{4,(3,1),(2,2),(2,11),(21,1),(111,1),1111\}.$$
Thus $E^4,E^{3,1},\do,E^{1111}$ in $\Irr(\WW)$, $E^4_q,E^{3,1}_q,\do,E^{1111}_q$ in $\Irr(\WW)$, (representations
of $\ch_q$) and $\r_4,\r_{3,1},\do,\r_{1111}$ (representations of $G(F_q)$) are defined.

Following \cite{\MA}, we index the irreducible representations of $\WW_0$ (of type $G_2$) as 
$$1_0,2_1,2_2,1'_3,1''_3,1_6.$$
This list is taken as the set $J$ so that $\fe_{1_0},\do,fe_{1_6}$ are in $\Irr(\WW_0)$; for $j\in J$ we write 
$R_j$ instead of $R_{\fe_j}$. From \cite{\MA, Theorem 6} we see that the following equalities (referred to as 
$\sp$) hold in $\cg$.    
$$\r_4=R_{1_0},$$
$$\r_{3,1}=R_{1'_3},\r_{111,1}=R_{1''_3},$$
$$\r_{2,2}=(1/2)R_{2_1}+(1/2)R_{2_2}+\Xi_1,$$
$$\r_{2,11}=(1/2)R_{2_1}-(1/2)R_{2_2}+\Xi_2,$$
$$\r_{21,1}=(1/2)R_{2_1}+(1/2)R_{2_2}+\Xi_3,$$
$$\r_{1111}=R_{1_6},$$
where $\Xi_1,\Xi_2,\Xi_3\in\cg_0$. Note that our $\r_i$ are defined over $\QQ$ hence they are the same as the 
extensions considered in \cite{\MA}.

In \cite{\MA, Theorem 10}, the values $R_j(g)$ are explicitly computed (as polynomials in $q$) for any $j\in J$
and any unipotent element in $D(F_q)$. The analogue of \cite{\MAA, Proposition 7} holds here with a similar 
proof; it implies that $\sum_{g\in\g(F_q)}\Xi(g)=0$ for any $\g\in\uD$, $\Xi\in\cg_0$. Hence we see from $\sp$ 
that

(a) {\it $\sum_{g\in\g(F_q)}\r_i(g)$ is explicitly known (as a polynomial in $q$) for any $\g\in\uD$ and any 
$i\in I_{ex}$.}
\nl
From 2.1(b) we have for any $\g\in\uD$ and any $w\in\WW_{D-min}$:
$$|\fB_w^\g(F_q)|=\sum_{g\in\g(F_q)}
\tr(S_gT_w:\cf@>>>\cf)=\sum_{i\in I_{ex}}(\sum_{g\in\g(F_q)}\tr(g,\r_i))\tr(V_iT_w,E^i_q).$$
Here the last sum can be calculated as an explicit polynomial in $q$ using (a) and the known values of the
polynomials $\tr(V_iT_w,E^i_q)$ (available through the CHEVIE package). This calculation was performed using a 
computer. The explicit knowledge of the quantities $|\{(g,B)\in\g(F_q)\T\cb(F_q);(B,gBg\i)\in\co_w\}|$ as 
polynomials in $q$ allows us to decide when $\g\in\Si_{w,D}$ (this holds precisely when the polynomial is 
nonzero). We can then verify that Theorem 1.3 holds in our case.

The calculation of the sum above yields in particular:
$$|\fB_w^\g(F_q)|=|G^0(F_q)|$$
whenever $C\in\uWW^{el}$, $w\in C_{min}$ and $\g=\Ph(C)$. From this and from \cite{\XW, 0.3(b)} we deduce as in
2.3 that 

(b) {\it $\fB_w^\g$ is a single $G^0$-orbit if $C\in\uWW^{el}$, $w\in C_{min}$ and $\g=\Ph(C)$.}
\nl
We now describe explicitly the map $\Ph$ in our case. We parametrize $\uWW_D$ as in the CHEVIE package. Thus we 
view $\WW$ as a subgroup of a Weyl group $W'$ of type $F_4$ which is normalized by an element $\o$ of order $3$ 
in $W'$ such that $\o y\o\i=\e_D(y)$ for $y\in\WW$. If $C\in\uWW_D$ then $C\o$ is contained in a unique conjugacy
class $C'$ of $W'$ and we give $C$ the same name as that given in \cite{\CA} to $C'$ (we also add an ${}^!$ to 
the name of an $\e_D$-elliptic conjugacy class). The objects of $\uD$ are denoted as 
$\g_{14},\g_8,\g_6,\g_4,\g_2$. Here the subscript denotes the codimension of the class (it is taken from 
\cite{\MA, Table VIII}). The map $\Ph$ is as follows:

$\tA_2\m\g_{14}$

$\tA_2A_2^!,\tA_2A_1\m\g_8$

$C_3A_1^!, C_3\m\g_6$

$F_4(a_1)^!\m\g_4\text{ dist}$   

$F_4^!\m\g_2\text{ dist}$.
\nl
Here we have indicated the distinguished unipotent classes by "dist".

\head 3. Even full orthogonal groups\endhead
\subhead 3.1\endsubhead
In this section we fix a $\kk$-vector space $V$ of finite even dimension $\nn=2n\ge4$ with a fixed nondegenerate
symplectic form $(,):V\T V@>>>\kk$ and a fixed quadratic form $Q:V@>>>\kk$ such that $(x,y)=Q(x+y)-Q(x)-Q(y)$ for
$x,y\in V$. For any subspace $V'$ of $V$ we set $V'{}^\pe=\{x\in V;(x,V')=0\}$. Let $\cl$ be the set of subspaces
$H$ of $V$ such that $\dim(H)=n$ and $Q|_H=0$. Let $Is(V)$ be the subgroup of $GL(V)$ consisting of all 
$g\in GL(V)$ such that $Q(gx)=Q(x)$ for all $x\in V$. For $j\in\{0,1\}$ let $G^j$ be the set of all vector space 
isomorphisms $g:V@>\si>>V$ such that for any $H\in\cl$ we have $\dim(g(H)\cap H)=n-j\mod2$. Note that 
$Is(V)=G^0\sqc G^1$ is an algebraic group with identity component $G^0$. In this section we assume that
$G=Is(V),D=G^1$. Then $G,D$ are as in 1.1(a).

\subhead 3.2\endsubhead
Let $\cf$ be the set of all sequences $V_*=(0=V_0\sub V_1\sub V_2\sub\do\sub V_\nn=V)$ of subspaces of $V$ such 
that $\dim V_i=i$ for $i\in[0,\nn]$, $Q|_{V_i}=0$ and $V_i^\pe=V_{\nn-i}$ for all $i\in[0,n]$. There is a unique 
involution $V_*\m\tV_*$ of $\cf$ where $\tV_i=V_i$ for $i\ne n$, $\tV_n\ne V_n$.

For $g\in Is(V)$, $V_*\in\cf$ we define $g\cdot V_*\in\cf$ by setting for any $i\in[0,\nn]$:
$(g\cdot V_*)_i=gV_i$. This defines an action of $Is(V)$ on $\cf$.

\subhead 3.3\endsubhead
Let $W$ be the group of permutations of $[1,\nn]$ which commute with the involution $i\m\nn-i+1$ of $[1,\nn]$. We
define a map $\cf\T\cf@>>>W$, $(V_*,V'_*)\m a_{V_*,V'_*}$ as in \cite{\WE, 1.4}; then $(V_*,V'_*)\m a_{V_*,V'_*}$
defines a bijection from the set of $Is(V)$-orbits on $\cf\T\cf$ (for the diagonal action) to $W$. For $w\in W$ 
let $\co_w$ be the $Is(V)$-orbit on $\cf\T\cf$ corresponding to $w$. Define $s_n\in W$ by $a_{V_*,\tV_*}=s_n$ for
any $V_*\in\cf$. We define $W'$ as the group of even permutations in $W$ (a subgroup of index $2$ of $W$). We 
view $W$ (resp. $W'$) as a Coxeter group of type $B_n$ (resp. $D_n$) as in \cite{\WE, 1.4}. We fix one of the two
$G^0$-orbits on $\cf$; we call it $\cf'$.

For any $V_*\in\cf$ we set $B_{V_*}=\{g\in G^0;g\cdot V_*=V_*\}$, a Borel subgroup of $G^0$. Then $V_*\m B_{V_*}$
is an isomorphism $\cf'@>\si>>\cb$, the variety of Borel subgroups of $G^0$. We identify $W'$ with $\WW$, the 
Weyl group of $G^0$, as in \cite{\WE, 1.5}. In our case we have $\e_D(w)=s_nws_n$ for $w\in\WW$. If $V_*\in\cf$ 
we have $B_{V_*}=B_{\tV_*}$.

\subhead 3.4\endsubhead
Until the end of 3.7 we fix $p_*=(p_1\ge p_2\ge\do\ge p_\s)$, a sequence of integers $\ge1$ such that $\s$ is odd
and $p_1+p_2+\do+p_\s=n$. Define $w_{p_*}\in W$ as in \cite{\WE, 1.6}. Note that $w_{p_*}\in W-W'$. Using 
\cite{\WE, 2.2(a)} and with notation of \cite{\WE, 1.4} we have
$$\align&s_nw_{p_*}\i=(s_{n-1})\do(s_{n-p_\s+1})\T\\&
(s_{n-p_\s}\do s_{n-1}\t_{n-1}\do s_{n-p_\s})(s_{n-p_\s-1})(s_{n-p_\s-2})\do(s_{n-p_\s-p_{\s-1}+1})\T\\&
(s_{n-p_\s-p_{\s-1}}\do s_{n-1}\ts_{n-1}\do s_{n-p_\s-p_{\s-1}})(s_{n-p_\s-p_{\s-1}-1})(s_{n-p_\s-p_{\s-1}-2})\\&
\do(s_{n-p_\s-p_{\s-1}-p_{\s-2}+1})\T\\&
\do\\&(s_{n-p_\s-\do-p_2}\do s_{n-1}\ts_{n-1}\do s_{n-p_\s-\do-p_2})(s_{n-p_\s-\do-p_2-1})\\&
(s_{n-p_\s-\do-p_2-2})\do(s_{n-p_\s-\do-p_1+1}).\endalign$$
Note that the length of $s_nw_{p_*}\i$ in $W'$ is equal to the length of $w_{p_*}\i$ in $W$ minus $\s$; hence it 
is $2\sum_{v=1}^{\s-1}vp_{v+1}+n-\s$. Thus $w_{p_*}s_n$ has minimal length in its $\e_D$-conjugacy class in 
$\WW=W'$ (see \cite{\GKP}). 

For any $g\in G^1$ let $X_g$ be the set of all $(V_*,V'_*)\in\cf\T\cf$ such that $a_{V_*,V'_*}=w_{p_*}$
(or equivalently $a_{V_*,\tV'_*}=w_{p_*}s_n$) and $g\cdot V_*=V'_*$. We have the following result \cite{\WE, 3.3}.

\proclaim{Proposition 3.5}Let $g\in G^1$ and let $(V_*,V'_*)\in X_g$. There exist vectors \lb $v_1,v_2,\do,v_\s$
in $V$ (each $v_i$ being unique up to multiplication by $\pm1$ such that for any $r\in[1,\s]$ the following hold.

(i) $V_{p_1+\do+p_{r-1}+i}=S(g^jv_k;k\in[1,r-1],j\in[0,p_k-1]\text{ or }k=r,j\in[0,i-1]\}$ for $i\in[0,p_r]$;

(ii) $(g^iv_t,v_r)=0$ for any $1\le t<r$, $i\in[-p_t,p_t-1]$;

(iii) $(v_r,g^iv_r)=0$ for $i\in[-p_r+1,p_r-1]$, $Q(v_r)=0$ and $(v_r,g^{p_r}v_r)=1$;

(iv) setting $E_r=S(g^{-p_t+i}v_t;t\in[1,r],i\in[0,p_t-1])$ we have $V=V_{p_{\le r}}\op E_r^\pe$.

(v) the vectors $(g^jv_t)_{t\in[1,\s],j\in[-p_t,p_t-1]}$ form a basis of $V$.
\endproclaim

\subhead 3.6\endsubhead
Let $g\in G^1$. Let $\cs_g$ be the set of all sequences $L_1,L_2,\do,L_\s$ of lines in $V$ such that for any 
$r\in[1,\s]$ we have

(i) $(g^iL_t,L_r)=0$ for any $1\le t<r$, $i\in[-p_t,p_t-1]$;

(ii) $(L_r,g^iL_r)=0$ for $i\in[-p_r+1,p_r-1]$, $Q(v_r)=0$;

(iii) $(L_r,g^{p_r}L_r)\ne0$.
\nl
We note:

(a) {\it if $L_1,L_2,\do,L_\s$ is in $\cs_g$ then the lines $\{g^{-p_t+i}L_t;t\in[1,\s],i\in[0,2p_t-1]\}$ form a 
direct sum decomposition of $V$}.
\nl
The proof is the same as that of 3.5(v).

Note that the assignment $(V_*,V'_*)\m(L_1,L_2,\do,L_\s)$  (where $L_i$ is spanned by $v_i$ as in 3.5) defines a 
bijection 

(b) $X_g@>\si>>\cs_g$.

\subhead 3.7\endsubhead
In the remainder of this section we assume that $p=2$ so that $G,D$ are as in 1.1(a),(b). Let $\tV=V\op\kk$. We 
define a quadratic form $\ti Q:\tV@>>>\kk$ by $\ti Q(x,e)=Q(x)+e^2$ where $x\in V,e\in\kk$. Let $(,)'$ be the 
symplectic form on $\tV$ attached to $\ti Q$. From \cite{\WE, 2.8}, \cite{\WE, 3.3}, \cite{\WE, 3.5(c)} we see 
that there exists a unipotent isometry $\tg:\tV@>>>\tV$ of $\ti Q$ with Jordan blocks of sizes 
$2p_1,2p_2,\do,2p_\s,1$ and vectors $v_1,v_2,\do,v_{\s+1}$ in $\tV$ so that the following hold.

$(\tg^iv_t,v_r)'=0$ for any $1\le t<r\le\s+1$, $i\in[-p_t,p_t-1]$;

$(v_r,\tg^iv_r)'=0$ for any $r\in[1,\s]$, $i\in[-p_r+1,p_r-1]$;

$\ti Q(v_r)=0$ and $(v_r,\tg^{p_r}v_r)'=1$ for any $r\in[1,\s]$;

$\ti Q(v_{\s+1})=1$;

for any $t\in[1,\s]$, the subspace spanned by $(\tg^{-p_t+i}v_t)_{i\in[0,2p_t-1]}$ is $\tg$-stable and 
$\tg v_{\s+1}=v_{\s+1}$;

the vectors $(\tg^jv_t)_{t\in[1,\s],j\in[-p_t,p_t-1]}$ together with $v_{\s+1}$ form a basis of $V$.
\nl
Let $V'$ be the subspace of $\tV$ spanned by $(\tg^jv_t)_{t\in[1,\s],j\in[-p_t,p_t-1]}$. Clearly, $(,)'$ is 
nondegenerate when restricted to $V'$ and $V'$ is a $\tg$-stable hyperplane in $\tV$ such that $\tg:V'@>>>V'$ is 
a unipotent isometry of $\ti Q|_{V'}$ with Jordan blocks of sizes $2p_1,2p_2,\do,2p_\s$.

Now the lines spanned by $v_1,v_2,\do,v_\s$ satisfy the definition of $\cs_{\tg|_{V'}}$ where $\cs_{\tg|_{V'}}$ 
is defined in terms of  $\tg|_{V'},V'$ in the same way as $\cs_g$ was defined in terms of $g,V$ in 3.6. Thus 
$\cs_{\tg|_{V'}}\ne\em$. Since $V',\tQ_{V'}$ is isomorphic to $V,Q$, it follows that there exists a unipotent 
isometry $g:V@>>>V$ of $Q$ with Jordan blocks of sizes $2p_1,2p_2,\do,2p_\s$ such that $\cs_g\ne\em$. Using 
3.6(b) we deduce that for this $g$ we have $X_g\ne\em$. Note that $g$ is necessarily in $G^1$.

\subhead 3.8\endsubhead
Let $(V_*,V'_*)\in\cf\T\cf$ be such that $a_{V_*,V'_*}=w_{p_*}$ or equivalently $a_{V_*,\tV_*}=w_{p_*}s_n$. We 
can assume that $V_*\in\cf'$. Then $\tV_*\in\cf'$.

Let $\g$ be a unipotent $G^0$-conjugacy class in $G^1$ such that $g\cdot V_*=V'_*$ for some $g\in\g$ and such 
that some/any $g\in\g$ has Jordan blocks of sizes $2p_1,2p_2,\do,2p_\s$. (Such $\g$ exists by 3.7.) Let $\g'$ be 
a unipotent $G^0$-conjugacy class in $G^1$ such that $g'\cdot V_*=V'_*$ for some $g'\in\g'$. We show:

(a) {\it $\g$ is contained in the closure of $\g'$.}
\nl
Let $Sp(V)$ the group of all automorphisms of $V$ that preserve $(,)$. We have $G^1\sub Sp(V)$. Let $\g_1,\g'_1$ 
be the $Sp(V)$-conjugacy class containing $\g,\g'$ respectively. By \cite{\LX, 1.3}, $\g_1$ is contained in the 
closure of $\g'_1$ in $Sp(V)$ and then, using \cite{\SP, II,8.2}, we see that $\g$ is contained in the closure of
$\g'$ in $G^1$. This proves (a).

From (a) we see that 1.3(a) holds for any $\e_D$-elliptic conjugacy class $C$ in $\WW$. We also see that 1.3(b),
1.3(d) hold for $G,D$ (we use the description of distinguished unipotent classes of $D$ given in \cite{\SP}).
This completes the verification of 1.4 for our $G,D$.

\head 4. Bilinear forms\endhead
\subhead 4.1\endsubhead
In this and the next two sections we assume that we are in one of the following two cases:

$\kk$ is an algebraically closed field and we set $q=1$, or

$\kk$ is an algebraic closure of a finite field $F_q$ with $q$ elements (we then have $q>1$).
\nl
Let $V$ be a $\kk$-vector space of dimension $n\ge1$; let $V^*$ be the dual vector space. For $x\in V,\x\in V^*$ 
we set $(x,\x)=\x(x)\in\kk$. For any subspace $V'$ of $V$ we set $V'{}^\pe=\{\x\in V^*;(V',\x)=0\}$; for any 
subspace $U$ of $V^*$ we set $U^\pe=\{x\in V;(x,U)=0\}$. For any $j\in\ZZ_{odd}$ let $G^j$ be the set of all 
group isomorphisms $g:V@>\si>>V^*$ such that $g(\l x)=\l^{q^j}g(x)$ for $x\in V,\l\in\kk$. For any $j\in\ZZ_{ev}$
let $G^j$ be the set of all group isomorphisms $g:V@>\si>>V$ such that $g(\l x)=\l^{q^j}g(x)$ for 
$x\in V,\l\in\kk$. We set $D=G^1$. For any $g\in G^j$ ($j$ odd) 
we define a group isomorphism $\chg:V^*@>\si>>V$ by $(\chg(\x),g(x))=(x,\x)^{q^j}$ for all $x\in V,\x\in V^*$; we
have $\chg(\l\x)=\l^{q^j}\chg(\x)$ for $\x\in V^*,\l\in\kk$. For any $g\in G^j$ ($j$ even) we define a group 
isomorphism $\chg:V^*@>\si>>V^*$ by $(g(x),\chg(\x))=(x,\x)^{q^j}$ for all $x\in V,\x\in V^*$; we have 
$\chg(\l\x)=\l^{q^j}\chg(\x)$ for $\x\in V^*,\l\in\kk$. 
For $g\in G^j$, $g'\in G^{j'}$ we define $g*g'\in G^{j+j'}$ by $g*g'=gg'$ if $j'$ is even; $g*g'=\chg g'$ if $j'$
is odd. If, in addition, $g''\in G^{j''}$, then $(g*g')*g''=g*(g'*g'')$ is equal to $gg'g''$ if $j',j''$ are 
even; to $\chg g'g''$ if $j'$ is odd, $j''$ is even; to $g\chg'g''$ if $j',j''$ are odd; to $\chg\chg'g''$ if 
$j'$ is even, $j''$ is odd. For $g\in G^j$ we define $g^{*(-1)}\in G^{-j}$ by $g^{*(-1)}=(\chg)\i$ if $j$ is odd 
and $g^{*(-1)}=g\i$ if $j$ is even. We have  $g^{*(-1)}*g=g*g^{*(-1)}=1\in G^0$. We see that if $q>1$, 
$\sqc_{j\in\ZZ}G^j$ is a group with multiplication $*$; we assume that it is the given group $G$ in 1.1, that
$G^0=GL(V)$ is the normal subgroup of $G$ given in 1.1 and $D$ is the given $G^0$-coset given in 1.1.
If $q=1$, $G^j$ depends only on the parity of $j$ and any element of $G^j$ is 
a vector space isomorphism; also, $*$ defines a group structure on $G^0\sqc G^1$ which becomes an algebraic group
with identity component $G^0=GL(V)$; we assume that $G^0\sqc G^1$ is the given group $G$ in 1.1, that
$G^0=GL(V)$ is the normal subgroup of $G$ given in 1.1 and $D$ is the given $G^0$-coset given in 1.1 (in this case
$G,D$ are as in 1.1(a)).

If $g\in G^1$ then the $j$-th power of $g$ in $G$ is denoted by $g^{*j}$; we have $g^{*j}\in G^j$; this agrees 
with the earlier definition of $g^{*(-1)}$. From the definitions we have $g^{*a}g^{*b}=g^{*(a+b)}$ if $b$ is even.
Also,
$$(x,gx')=(x',g^{*(-1)}x)^q\text{ for }x,x'\in V\tag a$$.
This implies
$$(g^{*(a+c)}x,g^{*(b+c)}x')=(g^{*a}x,g^{*b}x')^{q^c}$$
for $a,c$ even, $b$ odd and $x,x'\in V$.

\subhead 4.2\endsubhead
Let $\cf$ be the set of sequences $V_*=(0=V_0\sub V_1\sub V_2\sub\do\sub V_n=V)$ of subspaces of $V$ such that 
$\dim V_i=i$ for $i\in[0,n]$. We identify $\cf$ with $\cb$ by
$V_*\m B_{V_*}=\{x\in GL(V);xV_i=V_i\text{ for all }i\}$.

For $g\in G^j$, $V_*\in\cf$ we define $g\cdot V_*\in\cf$ by seeting for any $i\in[0,n]$:
$(g\cdot V_*)_i=gV_i$ (if $j$ is even), $(g\cdot V_*)_i=(gV_{n-i})^\pe$ (if $j$ is odd). If $g\in G^j$, 
$g'\in G^{j'}$ and $V_*\in\cf$, then $g\cdot(g'\cdot V_*)=(g*g')\cdot V_*$. (We use that $(gW)^\pe=\chg(W^\pe)$ 
for any subspace $W$ of $V$.) Thus, $g,V_*\m g\cdot V_*$ is a $G$-action on $\cf$. Under the identification
$\cf=\cb$, this becomes the $G$-action $g:B\m g*B*g^{*(-1)}$ on $\cb$.

\subhead 4.3\endsubhead
Let $S_n$ be the group of permutations of $[1,n]$. For any permutation $w\in S_n$ let $\co_w$ be the set of pairs
$(V_*,V'_*)\in\cf\T\cf$ such that $\dim(V'_j\cap V_{w(j)})=\dim(V'_j\cap V_{w(j)-1})+1$ for all $j\in[1,n]$. Note
that if $(V_*,V'_*)\in\co_w$ and $j\in[1,n]$ then
$$V_{w(j)}=(V'_j\cap V_{w(j)})+V_{w(j)-1}.\tag a$$
Now $\co_w$ is a single $G^0$-orbit on $\cf\T\cf$ (for the diagonal action). Hence $w$ can be viewed as an 
element of $\WW$. This identifies $\WW$ with $S_n$.

Define $\ww\in S_n$ by $j\m n+1-j$. From the definitions we see that

(b) {\it if $(V_*,V'_*)\in\co_w$ then $(V'_*,V_*)\in\co_{w\i}$.}

(c) {\it if $(V_*,V'_*)\in\co_w$ and $h\in G^1$ then $(h\cdot V_*,h\cdot V'_*)\in\co_{\ww w\ww}$.}
\nl
For $w\in\WW$ we have $\e_D(w)=\ww w\ww$. 

\subhead 4.4\endsubhead
Let $p_*=(p_1\ge p_2\ge\do\ge p_\s)$ be a descending sequence in $\ZZ_{>0}$ such that $2(p_1+p_2+\do+p_\s)=n+\s$.
Let $z_{p_*}$ be the permutation of $[1,n]$ such that for any $r\in[1,\s]$ we have
$$(p_1-1)+\do+(p_{r-1}-1)+i\m p_1+\do+p_{r-1}+i+1\qua (i\in[1,p_r-1]);$$
$$n-(p_1+\do+p_r-1)\m p_1+\do+p_{r-1}+1;$$
$$n-((p_1+\do+p_{r-1})+i)\m n-((p_1-1)+\do+(p_{r-1}-1)+i)\qua (i\in[0,p_r-2]).$$
From the results in \cite{\GKP} we see that $z_{p_*}$ has minimal legth in its $\e_D$-conjugacy class in 
$\WW=W'$. In the remainder of this section we write $w$ instead of $z_{p_*}$. 

The composition $\ww w$ is the permutation of $[1,n]$ such that for any $r\in[1,\s]$ we have
$$(p_1-1)+\do+(p_{r-1}-1)+i\m n-(p_1+\do+p_{r-1}+i)\qua (i\in[1,p_r-1]);$$
$$n-(p_1+\do+p_{r-1}+p_r-1)\m n-(p_1+\do+p_{r-1});$$
$$n-((p_1+\do+p_{r-1})+i)\m(p_1-1)+\do+(p_{r-1}-1)+i+1\qua (i\in[0,p_r-2]).$$
This is the product of $\s$ disjoint cycles of size $2p_1-1,2p_2-1,\do,2p_\s-1$:
$$1\m n-1\m 2\m n-2\m 3\m n-3\m\do\m p_1-1\m n-(p_1-1)\m n\m 1,$$
$$\align&p_1\m n-(p_1+1)\m p_1+1\m n-(p_1+2)\m\do\m p_1+p_2-2\m \\&n-(p_1+p_2-1)\m n-p_1\m p_1,\do.\endalign$$
The permutation $\ww w\i\ww$ is given by
$$(p_1-1)+\do+(p_{r-1}-1)+i\m(p_1+\do+p_{r-1})+i\qua (i\in[1,p_r-1]),$$
$$n-(p_1+\do+p_{r-1})\m p_1+\do+p_r,$$
$$n-(p_1+\do+p_{r-1}+i)\m n+1-((p_1-1)+\do+(p_{r-1}-1)+i)\qua (i\in[1,p_r-1]).$$
The following equality is a special case of 4.3(a):
$$\align&V_{p_1+\do+p_{r-1}+i+1}\\&=(V'_{(p_1-1)+\do+(p_{r-1}-1)+i}\cap V_{p_1+\do+p_{r-1}+i+1})
+V_{p_1+\do+p_{r-1}+i}\tag a\endalign$$
if $r\in[1,\s],i\in[0,p_r-1]$.

\subhead 4.5\endsubhead
Let $(V_*,V'_*)\in\co_w$. Using the equality 
$$\dim(V_j\cap V'_k)=\sha(h\in[1,k];w(h)\le j\}\tag a$$
we see that for any $r\in[1,\s]$ we have
$$\dim(V'_{(p_1-1)+\do+(p_{r-1}-1)+i}\cap V_{p_1+\do+p_{r-1}+i+1})=(p_1-1)+\do+(p_{r-1}-1)+i,$$
$$\dim(V'_{(p_1-1)+\do+(p_{r-1}-1)+i}\cap V_{p_1+\do+p_{r-1}+i}=(p_1-1)+\do+(p_{r-1}-1)+i-1$$
for $i\in[1,p_r-1]$;
$$\dim(V'_{n-((p_1+\do+p_{r-1})+i)}\cap V_{n-((p_1-1)+\do+(p_{r-1}-1)+i)})=n-((p_1+\do+p_{r-1})+i),$$
$$\align&\dim(V'_{n-((p_1+\do+p_{r-1})+i)}\cap V_{n-((p_1-1)+\do+(p_{r-1}-1)+i)-1})\\&=n-((p_1+\do+p_{r-1})+i)-1
\endalign$$
for $i\in[0,p_r-2]$;
$$\dim(V'_{n-(p_1+\do+p_r-1)}\cap V_{p_1+\do+p_{r-1}+1})=(p_1-1)+\do+(p_{r-1}-1)+1,$$
$$\dim(V'_{n-(p_1+\do+p_r-1)}\cap V_{p_1+\do+p_{r-1}})=(p_1-1)+\do+(p_{r-1}-1);$$
in particular we have
$$V'_{(p_1-1)+\do+(p_{r-1}-1)+i}\sub V_{p_1+\do+p_{r-1}+i+1}\tag b$$
for $i\in[1,p_r-1]$.
Similarly for $(V_*,V''_*)\in\co_{\ww w\i\ww}$ and any $r\in[1,\s],i\in[1,p_r-1]$ we have
$$\dim(V''_{(p_1-1)+\do+(p_{r-1}-1)+i}\cap V_{p_1+\do+p_{r-1}+i})=(p_1-1)+\do+(p_{r-1}-1)+i;$$
hence 
$$V''_{(p_1-1)+\do+(p_{r-1}-1)+i}\sub V_{p_1+\do+p_{r-1}+i}.\tag c$$

\subhead 4.6\endsubhead
Let $g\in G^1$. Let 
$$X_g=\{(V_*,V'_*)\in\in\co_w;V'_*=g\cdot V_*\},$$ 
$$X'_{g\i}=\{(V_*,V''_*)\in\co_{\ww w\i\ww};V''_*=g\i\cdot V_*\}.$$
We show:

(a) {\it If $(V_*,V'_*)\in X_g$ and $V''_*=g\i\cdot V_*$ then $(V_*,V''_*)\in X'_{g\i}$.}
\nl
Indeed, from $(V_*,g\cdot V_*)\in\co_w$ we deduce using 4.3(c) that $(g\i\cdot V_*,V_*)\in\co_{\ww w\ww}$ which
implies (a) in view of 4.3(b).

In the following result we fix $(V_*,V'_*)\in X_g$.

\proclaim{Proposition 4.7}Let $u\in[1,\s]$. There exist vectors $v_1,v_2,\do,v_u$ in $V$ (each $v_i$ being unique
up to multiplication by an element in $\{\l\in\kk^*;\l^{q^{2p_i-1}+1}=1\}$) such that for any $r\in[1,u]$ the 
following hold.

(i) for $i\in[1,p_r]$ we have
$$\align&V_{p_1+\do+p_{r-1}+i}\\&=S(g^{*j}v_k; k\in[1,r-1],j\in[0,2p_k-2]_{ev}\text{ or }k=r,j\in[0,2i-1]_{ev})
\endalign$$ 
(in particular, the vectors $g^{*j}v_k (k\in[1,r],j\in[0,2p_k-2]_{ev})$ are linearly independent);

(ii) for $i\in[1,p_r-1]$ we have
$$\align&V'_{(p_1-1)+\do+(p_{r-1}-1)+i}\\&=S(g^{*j}v_k; 
k\in[1,r-1],j\in[1,2p_k-1]_{ev}\text{ or }k=r,j\in[1,2i]_{ev});\endalign$$ 

(iii) setting $E_r=S(g^{*j}v_k;k\in[1,r],j\in[-2p_k+1,-1]_{odd})\sub V^*$  we have $V=V_{p_1+\do+p_r}\op E_r^\pe$;

(iv) $(v_r,g^{*j}v_t)=0$ if $t\in[1,\s]$, $t<r$, $j\in[-2p_t+1,2p_t-2]_{odd}$;

(v) $(v_r,g^{*j}v_r)=0$ if $j\in[-2p_r+2,2p_r-2]_{odd}$;

(vi) $(v_r,g^{*(2p_r-1)v_r})\ra=1$.    
\endproclaim
We can assume that the result holds when $u$ is replaced by a strictly smaller number in $[1,u]$. (This 
assumption is empty when $u=1$.) In particular $v_1,\do,v_{u-1}$ are defined. By assumption we have 
$V=V_{p_1+\do+p_{u-1}}\op E_{u-1}^\pe$ hence \lb 
$V_{p_1+\do+p_{u-1}+1}\cap E_{u-1}^\pe$ is a line. (We set $E_0=0$ 
so that $E_0^\pe=V$.) Let $v_u$ be a nonzero vector on this line. 

We show that (iv) holds for $r=u$. From the induction hypothesis we have 
$$V'_{(p_1-1)+\do+(p_{u-1}-1)}=S(g^{*j}v_k; k\in[1,u-1],j\in[1,2p_k-1]_{ev}).$$ 
Hence 
$$(gV_{n-((p_1-1)+\do+(p_{u-1}-1)})^\pe\sub S(g^{*j}v_k; k\in[1,u-1],j\in[1,2p_k-1]_{ev}).$$ 
We have $V_{p_1+\do+p_{u-1}+1}\sub V_{n-((p_1-1)+\do+(p_{u-1}-1)}$ since 
$$p_1+\do+p_{u-1}+1\le n-((p_1-1)+\do+(p_{u-1}-1)).$$
Hence 
$$\align&gv_u\in gV_{p_1+\do+p_{u-1}+1}\sub gV_{n-((p_1-1)+\do+(p_{u-1}-1)}\\&=
(V'_{(p_1-1)+\do+(p_{u-1}-1)})^\pe=S(g^{*j}v_k; k\in[1,u-1],j\in[1,2p_k-1]_{ev})^\pe.\endalign$$
(The last equality uses the induction hypothesis.) Thus if $k\in[1,u-1], j\in[1,2p_k-1]_{ev}$ we have
$(g^{*j}v_k,gv_u)=0$ so that if $k\in[u-1],j\in[0,2p_k-2]_{odd}$ then $(v_u,g^{*j}v_k)=0$. Since 
$v_u\in E_{u-1}^\pe$ we have $(v_u,g^{*j}v_k)=0$ if $k\in[1,u-1],j\in[-2p_k+1,-1]_{odd}$. Thus 
$(v_u,g^{*j}v_k)=0$ if $k\in[1,u-1],j\in[-2p_k+1,2p_k-2]_{odd}$. This proves (iv).

We show that (i) and (ii) hold for $r=u$.
It is enough to show (i) when $i\in[1,p_u]$ and (ii) when $i\in[1,p_u-1]$. From the definition we have 
$$\align&V_{p_1+\do+p_{u-1}+1}=V_{p_1+\do+p_{u-1}}+\kk v_u\\&
=S(g^{*j}v_k;k\in[1,u-1],j\in[0,2p_k-2]_{ev}\text{ or }k=u,j=0).\endalign$$
We assume that for some $a\in[1,p_u-1]$, (i) is known for $i=1,2,\do,a$ and (ii) is known for $i=1,2,\do,a-1$.
It is enough to show that (i) holds $i=a+1$ and (ii) holds for $i=a$. By assumption we have 
$$V_{p_1+\do+p_{u-1}+a}=S(g^{*j}v_k; k\in[1,u-1],j\in[0,2p_k-2]_{ev}\text{ or }k=u,j\in[0,2a-1]_{ev}),$$
$$\align&V'_{(p_1-1)+\do+(p_{u-1}-1)+a-1}\\&
=S(g^{*j}v_k; k\in[1,u-1],j\in[1,2p_k-1]_{ev}\text{ or }k=u,j\in[1,2a-2]_{ev}).\endalign$$
We have $V_{p_1+\do+p_{u-1}}\sub V_{n-((p_1-1)+\do+(p_{u-1}-1)+a)}$ since 
$$p_1+\do+p_{u-1}\le n-((p_1-1)+\do+(p_{u-1}-1)+a).$$
Hence for $k\in[1,u-1]$ we have
$$\align&gv_k\in gV_{p_1+\do+p_{u-1}}\sub gV_{n-((p_1-1)+\do+(p_{u-1}-1)+a)}\\&
=(V'_{(p_1-1)+\do+(p_{u-1}-1)+a})^\pe.\endalign$$
We have 
$$\align&V'_{(p_1-1)+\do+(p_{u-1}-1)+a}\\&
=S(g^{*j}v_k; k\in[1,u-1],j\in[1,2p_k-1]_{ev}\text{ or }k=u,j\in[1,2a-2]_{ev})\op\kk\et\endalign$$
for some $\et\in V$ and by the previous sentence we have $(\et,gv_k)$ for $k\in[1,u-1]$. Since 
$$\align&V'_{(p_1-1)+\do+(p_{u-1}-1)+a}\sub g^{*2}V_{p_1+\do+p_{u-1}+a}\\&
=S(g^{*j}v_k; k\in[1,u-1],j\in[2,2p_k]_{ev}\text{ or }k=u,j\in[2,2a]_{ev}),\endalign$$
(see 4.6(a), 4.5(c)) we see that we can assume that $\et=\sum_{k\in[1,u-1]}c_kg^{*(2p_k)}v_k+c_ug^{*(2a)}v_u$ 
where $c_1,\do,c_u\in\kk$ are not all zero. Assume that $c_r\ne0$ for some $r\in[1,u-1]$; let $r_0$ be the 
smallest such $r$. We have 
$$\align&0=(\et,gv_{r_0})=c_{r_0}\la v_{r_0},g^{*(2p_{r_0}-1)}v_{r_0}\ra^q\\&+\sum_{k\in[r_0+1,u-1]}
c_k(v_k,g^{*(-2p_k+1)}v_{r_0})^{q^{2p_k}}+c_u\la v_u,g^{*(-2a+1)}v_{r_0}\ra^{q^{2a}}=c_{r_0}.\endalign$$
(We have used (iv) and that for $k>r_0$ we have $-2p_k+1\in[-2p_{r_0}+1,-1]$.) Thus $c_{r_0}=0$, a
contradiction. We see that $\et$ is a nonzero multiple of $g^{*(2a)}v_u$. We can assume that $\et=g^{*(2a)}v_u$; 
we see that (ii) holds for $i=a$.

We have $V_{p_1+\do+p_{u-1}+a+1}=V'_{(p_1-1)+\do+(p_{u-1}-1)+a}+V_{p_1+\do+p_{u-1}+a}$. (We use 4.4(a) and the 
inclusion $V'_{(p_1-1)+\do+(p_{u-1}-1)+a}\sub V_{p_1+\do+p_{u-1}+a+1}$ in 4.5(b).) Using this, (ii) for $i=a$ and
the induction hupothesis we see that 
$$\align&V_{p_1+\do+p_{u-1}+a+1}=S(g^{*j}v_k; k\in[1,u-1],j\in[1,2p_k-1]_{ev}\text{ or }k=u,j\in[1,2a]_{ev})\\&
+S(g^{*j}v_k; k\in[1,r-1],j\in[0,2p_k-2]_{ev}\text{ or }k=r,j\in[0,2a-1]_{ev})\\&=
S(g^{*j}v_k; k\in[1,u-1],j\in[0,2p_k-2]_{ev}\text{ or }k=u,j\in[0,2a+1]_{ev})\endalign$$
so that (i) holds for $i=a+1$. This proves (i) and (ii).

We show that (v) holds for $r=u$. 
We have $V_{p_1+\do+p_u}\sub V_{n-((p_1-1)+\do+(p_u-1))}$ since $p_1+\do+p_u\le n-((p_1-1)+\do+(p_u-1))$. Hence
$$gV_{p_1+\do+p_u}\sub gV_{n-((p_1-1)+\do+(p_u-1))}=(V'_{(p_1-1)+\do+(p_u-1)})^\pe.$$
Using this and (i),(ii) we deduce
$$S(g^{*j}v_k; k\in[1,u],j\in[1,2p_k-1]_{odd})\sub S(g^{*j}v_k; k\in[1,u],j\in[1,2p_k-1]_{ev})^\pe.$$
In particular if $j\in[1,2p_u-1]_{ev}$, $j'\in[1,2p_u-1]_{odd}$ then $(g^{*j}v_u,g^{*j'}v_u)=0$ and (v) follows.

We show that (vi) holds for $r=u$. From 4.5(a) we have 
$$\dim(V'_{n-(p_1+\do+p_u)}\cap V_{p_1+\do+p_{u-1}+1})=(p_1-1)+\do+(p_{u-1}-1).$$
Hence 
$$\dim(V'_{n-(p_1+\do+p_u)}+V_{p_1+\do+p_{u-1}+1})=n-(p_1-1)+\do+(p_u-1)$$ 
so that
$$\dim((V'_{n-(p_1+\do+p_u)})^\pe\cap V_{p_1+\do+p_{u-1}+1}^\pe)=(p_1-1)+\do+(p_u-1)$$
that is
$$\dim(gV_{p_1+\do+p_u}\cap V_{p_1+\do+p_{u-1}+1}^\pe)=(p_1-1)+\do+(p_u-1).$$
By (i) we have $S(g^{*j}v_k;k\in[1,u],j\in[1,2p_k-3]_{odd})\sub gV_{p_1+\do+p_u}$. By (i) and (iv),(v), we have 
$S(g^{*j}v_k;k\in[1,u],j\in[1,2p_k-3]_{odd})\sub V_{p_1+\do+p_{u-1}+1}^\pe$.
Hence 
$$S(g^{*j}v_k;k\in[1,u],j\in[1,2p_k-3]_{odd})\sub gV_{p_1+\do+p_u}\cap V_{p_1+\do+p_{u-1}+1}^\pe.$$
By (i) we have 
$$\dim S(g^{*j}v_k;k\in[1,u],j\in[1,2p_k-3]_{odd})=(p_1-1)+\do+(p_u-1).$$ 
Hence we must have
$$S(g^{*j}v_k;k\in[1,u],j\in[1,2p_k-3]_{odd})=gV_{p_1+\do+p_u}\cap V_{p_1+\do+p_{u-1}+1}^\pe.$$
By (i) we have 
$$\align&\dim S(g^{*j}v_k;k\in[1,u],j\in[1,2p_k-3]_{odd}\text{ or } k=u,j=2p_u-1)\\&=(p_1-1)+\do+(p_u-1)+1
\endalign$$
hence 
$$g^{*(2p_u-1)}v_u\n S(g^{*j}v_k;k\in[1,u],j\in[1,2p_k-3]_{odd})$$
so that
$$g^{*(2p_u-1)}v_u\n gV_{p_1+\do+p_u}\cap V_{p_1+\do+p_{u-1}+1}^\pe.$$
Since
$g^{*(2p_u-1)}v_u\in gV_{p_1+\do+p_u}$ (by (i)) it follows that $g^{*(2p_u-1)}v_u\n V_{p_1+\do+p_{u-1}+1}^\pe$. 
Using (i) we deduce
$$g^{*(2p_u-1)}v_u\n S(g^{*j}v_k;k\in[1,u-1],j\in[0,2p_k-2]_{ev}\text{ or }k=u,j=0)^\pe.$$
Since
$$g^{*(2p_u-1)}v_u\n S(g^{*j}v_k;k\in[1,u-1],j\in[0,2p_k-2]_{ev})^\pe$$
(see (iv)) we deduce that $(v_u,g^{*(2p_u-1)}v_u)\ne0$. Replacing $v_u$ by a scalar multiple we see that we can 
assume that $(v_u,g^{*(2p_u-1)}v_u)=1$. This proves (vi).

We show that (iii) holds for $r=u$. Note that $\dim(E_u)\le p_1+\do+p_r$; hence
$\dim(E_u^\pe)+\dim V_{p_1+\do+p_u}\ge n$. Using this we see that it is enough to show that
$E_u^\pe\cap V_{p_1+\do+p_u}=0$. Taking (i) into account we see that it is enough to verify the following 
statement:

($*$) {\it Let $f=\sum_{k\in[1,u],i\in[1,p_k]}c_{k,i}g^{*(2p_k-2i)}v_k$ ($c_{k,i}\in\kk$) be such that 
$(f,g^{*j'}v_{k'})=0$ for any $k'\in[1,u],j'\in[-2p_{k'}+1,-1]_{odd}$. Then $f=0$.}
\nl
Assume that not all $c_{k,i}$ are zero. Let $i_0=\min\{i;c_{k,i}\ne0\text{ for some }k\in[1,u]\}$. Let 
$X'=\{k\in[1,u];c_{k,i_0}\ne0\}$. We have $X'\ne\em$ and 
$$f=\sum_{r\in X'}c_{r,i_0}g^{*(2p_r-2i_0)}v_r+\sum_{r\in[1,u]}\sum_{i\in[i_0+1,p_r]}c_{r,i}g^{*(2p_r-2i)}v_r.$$
Let $r_0$ be the smallest number in $X'$. We have 
$$\align&0=(f,g^{*(-2i_0+1)}v_{r_0})=\sum_{r\in X'}c_{r,i_0}(g^{*(2p_r-2i_0)}v_r,g^{*(-2i_0+1)}v_{r_0})
\\&+\sum_{r\in[1,u]}\sum_{i\in[i_0+1,p_r]}c_{r,i}(g^{*(2p_r-2i)}v_r,g^{*(-2i_0+1)}v_{r_0}).\endalign$$
If $r\in X'$, $r\ne r_0$, we have $(g^{*(2p_r-2i_0)}v_r,g^{*(-2i_0+1)}v_{r_0})=(v_r,g^{*(-2p_r+1)}v_{r_0})=0$ 
(we use (iv); note that $r\ge r_0$ hence $p_r\le p_{r_0}$). If $r\in[1,u]$ and $i\in[i_0+1,p_r]$, we have 
$(g^{*(2p_r-2i)}v_r,g^{*(-2i_0+1)}v_{r_0})=0$ (we use (iv),(v); note that if $r>r_0$ we have 
$2(-p_r+i-i_0)+1\in [-2p_{r_0}+1,-1]$; if $r\le r_0$ we have $2(p_r-i+i_0)-1\in[0,2p_r-3]$.) Thus we have 
$$0=c_{r_0,i_0}(g^{*(2p_{r_0}-2i_0)}v_{r_0},g^{*(-2i_0+1)}v_{r_0})=
c_{r_0,i_0}(v_{r_0},g^{*(-2p_{r_0}+1)}v_{r_0})=c_{r_0,i_0}.$$
We see that $c_{r_0,i_0}=0$, a contradiction. 
This proves ($*$) hence (iii).

This completes the proof of existence part of the proposition. The uniqueness part follows from the proof of 
existence. The proposition is proved.

\subhead 4.8\endsubhead
Let $g\in G^1$. Let $\cs_g$ be the set of all sequences $L_1,L_2,\do,L_\s$ of lines in $V$ such that for any 
$r\in[1,\s]$ we have

(i) $(L_r,g^{*j}L_t)=0$ if $t\in[1,\s]$, $t<r$, $j\in[-2p_t+1,2p_t-2]_{odd}$;

(ii) $(L_r,g^{*j}L_r)=0$ if $j\in[-2p_r+2,2p_r-2]_{odd}$;

(iii) $(L_r,g^{*(2p_r-1)}L_r)\ne0$.
\nl
We show:

(a) {\it if $L_1,L_2,\do,L_\s$ is in $\cs_g$ then the lines $\{g^{*(-2p_k+2h)}L_k;k\in[1,\s],h\in[0,2p_k-2]\}$
form a direct sum decomposition of $V$}.
\nl
For $r\in[1,\s]$ we can find $v_r\in L_r$ such that $(v_r,g^{*(2p_r-1)}v_r)=1$. Assume that 
$f=\sum_{k\in[1,\s]}\sum_{h\in[0,2p_k-2]}c_{k,h}g^{*(-2p_k+2h)}v_k$ is equal to $0$ where $c_{k,h}\in\kk$ are not
all zero. Let $h_0=\min\{h;c_{k,h}\ne0\text{for some }k\in[1,\s]\}$. Let $X=\{k\in[1,\s];c_{k,h_0}\ne0\}$. We 
have $X\ne\em$ and 
$$f=\sum_{k\in X}c_{k,h_0}g^{*(-2p_k+2h_0)}v_k+\sum_{k\in[1,\s]}\sum_{h\in[h_0+1,2p_k-2]}g^{*(-2p_k+2h)}v_k.$$
Let $k_0$ be the largest number in $X$. We have
$$\align&0=(f,g^{*(2h_0-1)}v_{k_0})=\sum_{k\in X}c_{k,h_0}(g^{*(-2p_k+2h_0)}v_k,g^{*(2h_0-1)}v_{k_0})\\&
+\sum_{k\in[1,\s]}\sum_{h\in[h_0+1,2p_k-2]}(g^{*(-2p_k+2h)}v_k,g^{*(2h_0-1)}v_{k_0}).\endalign$$
If $k\in X$, $k\ne k_0$ we have $(g^{*(-2p_k+2h_0)}v_k,g^{*(2h_0-1)}v_{k_0})=0$ (using (i) and $k<k_0$). If 
$k\in[1,\s]$ and $h\in[h_0+1,2p_k-2]$ we have $(g^{*(-2p_k+2h)}v_k,g^{*(2h_0-1)}v_{k_0})=0$ (we use (i),(ii); 
note that if $k<k_0$, we have $-2p_k+2h-2h_0+1\in[-2p_k+1,2p_k-2]$; if $k\ge k_0$ we have
$2p_k-2h+2h_0-1\in[-2p_{k_0}+2,2p_{k_0}-2]$). We see that
$$0=c_{k_0,h_0}(g^{*(-2p_{k_0}+2h_0)}v_{k_0},g^{*(2h_0-1)}v_{k_0}).$$
Using (iii) we deduce $0=c_{k_0,h_0}$, a contradiction.

\subhead 4.9\endsubhead
Let $g\in G^1$. For any $L_1,L_2,\do,L_\s$ in $\cs_g$ we define subspaces $V_j,V'_j$ of $V$ ($j\in[1,n-1])$ as 
follows:
$$V_{p_1+\do+p_{r-1}+i}=\sum_{k\in[1,r-1],j\in[0,2p_k-2]_{ev}\text{ or }k=r,j\in[0,2i-1]_{ev}}g^{*j}L_k,$$
$r\in[1,\s],i\in[1,p_r]$;
$$V_{n-((p_1-1)+\do+(p_{r-1}-1)+i)}=
(\sum_{k\in[1,r-1],j\in[0,2p_k-2]_{odd}\text{ or }k=r,j\in[0,2i-1]_{odd}}g^{*j}L_k)^\pe,$$
$r\in[1,\s],i\in[1,p_r-1]$;
$$V'_{(p_1-1)+\do+(p_{r-1}-1)+i}=
\sum_{k\in[1,r-1],j\in[1,2p_k-1]_{ev}\text{ or }k=r,j\in[1,2i]_{ev}}g^{*j}L_k,$$
$r\in[1,\s],i\in[1,p_r-1]$;
$$V'_{n-(p_1+\do+p_{r-1}+i)}=
(\sum_{k\in[1,r-1],j\in[1,2p_k-1]_{odd}\text{ or }k=r,j\in[1,2i]_{odd}}g^{*j}L_k)^\pe,$$ 
$r\in[1,\s],i\in[1,p_r]$.
Note that the sums above are direct, by 4.8(a). Note also that $V_{p_1+\do+p_\s}=V_{n-((p_1-1)+\do+(p_\s-1))}$ is
defined in two different ways; similarly, $V'_{(p_1-1)+\do+(p_\s-1)}=V'_{n-(p_1+\do+p_\s)}$ is defined in two 
different ways; these two definitions are compatible by the definition of $\cs_g$. We set 
$V_n=V'_n=V,V_0=V'_0=0$. We have $V_*=(V_j)\in\cf$, $V'_*=(V'_j)\in\cf$ and $V'_j=(gV_{n-j})^\pe$ for all 
$j\in[1,n]$. 

For $r\in[1,\s],i\in[1,p_r-1]$ we have
$$\align&V'_{(p_1-1)+\do+(p_{r-1}-1)+i}\cap V_{p_1+\do+p_{r-1}+i+1}\\&=
\op_{k\in[1,r-1],j\in[1,2p_k-1]_{ev}\text{ or }k=r,j\in[1,2i]_{ev}}g^{*j}L_k\\&
\cap \op_{k\in[1,r-1],j\in[0,2p_k-2]_{ev}\text{ or }k=r,j\in[0,2i+1]_{ev}}g^{*j}L_k \\&
=\op_{k\in[1,r-1],j\in[1,2p_k-2]_{ev}\text{ or }k=r,j\in[1,2i]_{ev}}g^{*j}L_k\endalign$$ 
and this has dimension $(p_1-1)+\do+(p_{r-1}-1)+i$;
$$\align&V'_{(p_1-1)+\do+(p_{r-1}-1)+i}\cap V_{p_1+\do+p_{r-1}+i}\\&=
\op_{k\in[1,r-1],j\in[1,2p_k-1]_{ev}\text{ or }k=r,j\in[1,2i]_{ev}}g^{*j}L_k \\&
\cap \op_{k\in[1,r-1],j\in[0,2p_k-2]_{ev}\text{ or }k=r,j\in[0,2i-1]_{ev}}g^{*j}L_k \\&
=\op_{k\in[1,r-1],j\in[1,2p_k-2]_{ev}\text{ or }k=r,j\in[1,2i-1]_{ev}}g^{*j}L_k\endalign$$ 
and this has dimension $(p_1-1)+\do+(p_{r-1}-1)+i-1$. For $r\in[1,\s],i\in[0,p_r-2]$ we have
$$\align&V'_{n-(p_1+\do+p_{r-1}+i)}\cap V_{n-((p_1-1)+\do+(p_{r-1}-1)+i)}\\&
=(\op_{k\in[1,r-1],j\in[1,2p_k-1]_{odd}\text{ or }k=r,j\in[1,2i]_{odd}}g^{*j}L_k)^\pe\\&
\cap(\op_{k\in[1,r-1],j\in[0,2p_k-2]_{odd}\text{ or }k=r,j\in[0,2i-1]_{odd}}g^{*j}L_k)^\pe\\&
=(\op_{k\in[1,r-1],j\in[0,2p_k-1]_{odd}\text{ or }k=r,j\in[0,2i]_{odd}}g^{*j}L_k)^\pe\endalign$$
and this has dimension $n-(p_1+\do+p_{r-1}+i)$;
$$\align&V'_{n-(p_1+\do+p_{r-1}+i)}\cap V_{n-((p_1-1)+\do+(p_{r-1}-1)+i+1)}\\&
=(\op_{k\in[1,r-1],j\in[1,2p_k-1]_{odd}\text{ or }k=r,j\in[1,2i]_{odd}}g^{*j}L_k)^\pe\\&
\cap (\op_{k\in[1,r-1],j\in[0,2p_k-2]_{odd}\text{ or }k=r,j\in[0,2i+1]_{odd}}g^{*j}L_k)^\pe\\&
=\op_{k\in[1,r-1],j\in[0,2p_k-1]_{odd}\text{ or }k=r,j\in[0,2i+1]_{odd}}g^{*j}L_k)^\pe\endalign$$
and this has dimension $n-(p_1+\do+p_{r-1}+i+1)$. For $r\in[1,\s]$ we have
$$\align&V'_{n-(p_1+\do+p_r-1)}\cap V_{p_1+\do+p_{r-1}}\\&
=\op_{k\in[1,r-1],j\in[1,2p_k-1]_{odd}\text{ or }k=r,j\in[1,2p_r-2)]_{odd}}g^{*j}L_k)^\pe\\&
\cap \op_{k\in[1,r-1],j\in[0,2p_k-2]_{ev}}g^{*j}L_k\\&=
\op_{k\in[1,r-1],j\in[1,2p_k-2]_{ev}}g^{*j}L_k\endalign$$ 
and this has dimension $(p_1-1)+\do+(p_{r-1}-1)$;
$$\align&V'_{n-(p_1+\do+p_r-1)}\cap V_{p_1+\do+p_{r-1}+1}\\&
=(\op_{k\in[1,r-1],j\in[1,2p_k-1]_{odd}\text{ or }k=r,j\in[1,2p_r-2)]_{odd}}g^{*j}L_k)^\pe\\&
\cap\op_{k\in[1,r-1],j\in[0,2p_k-2]_{ev}\text{ or } k=r,j\in[0,1]_{ev}}g^{*j}L_k \\&
=\op_{k\in[1,r-1],j\in[1,2p_k-2]_{ev}\text{ or }k=r,j=0}g^{*j}L_k\endalign$$ 
and this has dimension $(p_1-1)+\do+(p_{r-1}-1)+1$. (We use the definition of $\cs_g$.) We see that 
$(V_*,V'_*)\in X_g$. Thus $(L_1,\do,L_\s)\m(V_*,V'_*)$ defines a morphism 
$$\cs_g@>>>X_g.\tag a$$
This is an isomorphism; the inverse is provided by $(V_*,V'_*)\m(L_1,\do,L_\s)$ where $L_r$ is the line spanned 
by $v_r$ as in 4.7 (with $u=\s$).

\head 5. Almost unipotent bilinear forms\endhead
\subhead 5.1\endsubhead
In this section we preserve the setup of 4.1 and we assume that $q=1$, that $G$ (which was fixed in 1.1) is as in
4.1 and $D=G^1$. Thus $G,D$ are as in 1.1(a). Note that $D$ does not contain unipotent elements of $G$, unless 
$p=2$. Let $D_{au}$ be the set of all $g\in D=G^1$ such that $g^{*2}:V@>>>V$ is unipotent (we then say that $g$ 
is {\it almost unipotent}). 

Until the end of 5.9 we fix $p_*=(p_1\ge p_2\ge\do\ge p_\s)$, a descending sequence in $\ZZ_{>0}$ such that 
$2(p_1+p_2+\do+p_\s)=n+\s$; note that the elements $z_{p_*}$ (see 4.4) with $p_*$ as above form a set of 
representatives for the elliptic $\e_D$-conjugacy classes in $\WW$. We fix a basis 
$\{z^t_i;t\in[1,\s],i\in[0,2p_t-2]\}$ of $V$. For any $r\in[1,\s],j\in[0,2p_r-2]$ we define $z'{}^r_j\in V^*$ by 
the following requirements:

$(z^t_i,z'{}^r_j)=0$ if $t\ne r$

$(z^r_i,z'{}^r_j)=0$ if $-p_r+1\le j-i\le p_r-2$

$(z^r_i,z'{}^r_j)=\bin{j-i+p_r-1}{j-i-p_r+1}$ if $j-i\ge p_r-1$,   

$(z^r_i,z'{}^r_j)=\bin{i-j+p_r-2}{i-j-p_r}$ if $j-i\le-p_r$.       
\nl
Clearly, the $n\T n$ matrix with entries $(z^t_i,z'{}^r_j)$ is nonsingular. Hence the $z'{}^r_j$ form a basis of 
$V^*$. Define $g:V@>\si>>V^*$ by $g(z^t_i)=z'{}^t_i$ for any $t\in[1,\s],i\in[0,2p_t-2]$. Define $g':V^*@>>>V$ by 

$g'(z'{}^r_j)=z^r_{j+1}$ for any $r\in[1,\s],j\in[0,2p_r-3]$,

$g'(z'{}^r_{2p_r-2})=\sum_{k\in[0,2p_r-2]}(-1)^k\bin{2p_r-1}{k}z^r_k$ for any $r\in[1,\s]$. 
\nl
We show that $g'=\chg$ that is, 

(a) $(g'(z'{}^r_j),g(z^t_i))=(z^t_i,z'{}^r_j)$ for all $r,j,t,i$. 
\nl
When $t\ne r$ both sides of (a) are zero. Thus we can assume that $r=t$. Assume first that $j<2p_t-2$. We must 
show that 

$(z^t_{j+1},z'{}^t_i)=(z^t_i,z'{}^t_j)$.
\nl
If $-p_t+1\le i-j-1\le p_t-2$ (or equivalently if $-p_t+1\le j-i\le p_t-2$) then the left hand side is $0$ and
the right hand side is $0$).

If $i-j-1\ge p_t-1$ (or equivalenty if $j-i\le-p_t$) then the left hand side is $\bin{i-j-1+p_t-1}{i-j-1-p_t+1}$ 
and the right hand side is $\bin{i-j+p_t-2}{i-j-p_t}$.

If $i-j-1\le-p_t$ (or equivalenty if $j-i\ge p_t-1$) then the left hand side is $\bin{j+1-i+p_t-2}{j+1-i-p_t}$ 
and the right hand side is $\bin{j-i+(p_t-1)}{j-i-(p_t-1)}$.

Next we assume that $j=2p_t-2$; we must check for any $i\in[0,2p_t-2]$ that
$$\sum_{k\in[0,2p_t-2]}(z^r_k,z'{}^t_i)(-1)^k\bin{2p_t-1}{k}=(z^t_i,z'{}^t_{2p_t-2})$$
that is,
$$\align&\sum_{k\in[0,2p_t-2];k\le i-(p_t-1)}(-1)^k\bin{i-k+(p_t-1)}{i-k-(p_t-1)}\bin{2p_t-1}{k}\\&
+\sum_{k\in[0,2p_t-2]; k\ge i+p_t}(-1)^k\bin{k-i+p_t-2}{k-i-p_t}\bin{2p_t-1}{k}\endalign$$
is equal to $0$ if $i\in[p_t,2p_t-1]$ and to $\bin{3p_t-3-i}{p_t-1-i}$ if $i\in[0,p_t-1]$. If $i=p_t-1$ the 
desired equality is $1=1$. It remains to show that
$$\sum_{k\in[i+p_t,2p_t-1]}(-1)^k\bin{k-i+p_t-2}{k-i-p_t}\bin{2p_t-1}{k}=0$$
if $i\in[0,p_t-2]$ and
$$\sum_{k\in[0,i-(p_t-1)]}(-1)^k\bin{i-k+(p_t-1)}{i-k-(p_t-1)}\bin{2p_t-1}{k}=0$$
if $i\in[p_t,2p_t-2]$.
Both of these equalities are special cases of the identity
$$\sum_{k,k'\ge0,k+k'=j}(-1)^k\bin{2p_t-2+k'}{k'}\bin{2p_t-1}{k}=0$$
for any $j\ge1$, which is easily verified. (We use the convention that $\bin{2p_t-1}{k}=0$ if $k>2p_t-1$.) This 
completes the proof of (a).

Now $g^{*2}:V@>>>V$ is given by  $z^t_i\m z^t_{i+1}$ for any $t\in[1,\s],i\in[0,2p_t-3]$ and 
$$z^t_{2p_t-2}\m\sum_{k\in[0,2p_t-2]}(-1)^k\bin{2p_t-1}{k}z^t_k$$
for any $t\in[1,\s]$. It follows that $g^{*2}:V@>>>V$ is unipotent (that is $g\in D_{au}$) with Jordan blocks of 
sizes $2p_1-1,2p_2-1,\do,2p_\s-1$.

For $r\in[1,\s]$ we set $v_r=g^{*2}z^r_{p_r-1}$. For $r\in[1,\s]$, $j\in[-2p_r+2,2p_r-2]_{odd}$ we have 
$0\le p_r-1+(j-1)/2\le 2p_r-3$ hence $z'{}^t_{p_r-1+(j-1)/2}$ is defined and
$$(v_r,g^{*j}v_r)=(g^{*2}z^r_{p_r-1},g^{*(j+2)}z^r_{p_r-1})=(z^r_{p_r-1},z'{}^t_{p_r-1+(j-1)/2})=0$$
(we use that $-p_r+1\le(j-1)/2\le p_r-2$). For $r\in[1,\s]$ we have
$$\align&(v_r,g^{*(2p_r-1)}v_r)=(g^{*2}z^r_{p_r-1},g^{*(2p_r-1+2)}z^r_{p_r-1})\\&
=(z^r_{p_r-1},g^{*(2p_r-1)}z^r_{p_r-1})=(z^r_{p_r-1},z'{}^r_{2p_r-2})=\bin{2p_r-2}{0}=1.\endalign$$
For $t,r\in[1,\s]$, $t\ne r$ and $j$ odd we have
$$(v_r,g^{*j}v_t)=(g^{*2}z^r_{p_r-1},g^{*(j+2)}z^t_{p_t-1})=(g^{*((j-1)/2)}z^r_{p_r-1},z'{}^t_{p_t-1})=0$$
since $g^{*((j-1)/2)}z^r_{p_r-1}$ is a linear combination of $z^r_0,z^r_1,\do,z^r_{2p_r-2}$. Thus, if $L_r$ is 
the line spanned by $v_r$,  we have $(L_1,L_2,\do,L_\s)\in\cs_g$ (see 4.8). Using the isomorphism 4.9(a) we see 
that
$$X_g\ne\em\tag a$$
(notation of 4.6).

\subhead 5.2\endsubhead
Let $g\in D_{au}$ be such that $X_g\ne\em$ (see 4.6); let $(V_*,V'_*)\in X_g$. Let $v_1,v_2,\do,v_\s$ be as in 
4.7 (with $u=\s$). For any $t\in[1,\s]$ let $W_t$ be the subspace of $V$ spanned by 
$(g^{*(-2p_r+2h)}v_r)_{r\in[0,t],h\in[0,2p_r-2]}$ and let $W'_t$ be the subspace of $V$ spanned by
$(g^{*(-2p_r+2h)}v_r)_{r\in[t+1,\s],h\in[0,2p_r-2]}$. From 4.8(a) we see that
$$V=W_t\op W'_t.$$
We set $N=g^{*2}-1:V@>>>V$, a nilpotent linear map. Let $n_1\ge n_2\ge\do\ge n_u$ be the sizes of the Jordan 
blocks of $N$; we set $n_i=0$ for $i>u$. For any $k\ge0$ we set
$$\L'_k=\sum_{r\in[1,\s]}\max(2p_r-1-k,0).$$
We show:

(a) {\it For any $k\ge1$ we have $\dim N^kV\ge\L'_k$.}
\nl
Note that $\dim N^kV=\sum_{i\ge1}\max(n_i-k,0)$. Applying \cite{\WE, 3.1(b)} with $x_1,\do,x_f$ given by 
$g^{*(-2p_1)}v_1,\do,g^{*(-2p_\s)}v_\s$ and using 4.8(a) we see that for any $c\ge1$ we have 
$(2p_1-1)+(2p_2-1)+\do+(2p_c-1)\le n_1+n_2+\do+n_c$. Hence for any $k\ge1$ we have 
$\sum_{r\in[1,\s]}\max(2p_r-1-k,0)\le\sum_{r\ge1}\max(n_r-k,0)$ and (a) follows.

We now assume that $k>0$ and $d\in[1,\s]$ is such that $2p_d-1\ge k$ and (if $d<\s$) $k\le2p_{d+1}-1$. Then
$\L'_k=\sum_{r\in[1,d]}(2p_r-1-k)$. We show:

(b) {\it If $\dim N^kV=\L'_k$ then $W_d,W'_d$ are $g^{*2}$-stable, $gW'_d=W_d^\pe$, $g^{*2}:W_d@>>>W_d$ has 
exactly $d$ Jordan blocks (each one has size $\ge k$) and $N^kW'_d=0$.}
\nl
For $r\in[1,\s]$ let $v'_r=g^{*(-2p_r)}v_r$; then $(g^{*(2h)}v'_r)_{h\in[0,2p_r-2]}$ is a basis of $X_r$ hence 
$(N^hv'_r)_{h\in[0,2p_r-2]}$ is a basis of $X_r$. For $r\in[1,d]$ let $Y_r$ be the subspace spanned by 
$N^hv'_r(h\in[k,2p_r-2])$. Note that $Y_r\sub N^kX_r$. Hence $\op_{r\in[1,d]}Y_r\sub N^kW_d\sub N^kV$. We have 
$\dim\op_{r\in[1,d]}Y_r=\sum_{r\in[1,d]}(2p_r-1-k)=\L'_k=\dim N^kV$. Hence $\op_{r\in[1,d]}Y_r=N^kW_d=N^kV$. We 
have $\op_{r\in[1,d]}Y_r\sub W_d$. Hence $N^kV\sub W_d$. We show that $NW_d\sub W_d$. Clearly $N$ maps the basis 
elements $N^hv'_r$ $(r\in[1,d],h\in[0,2p_r-3])$ into $W_d$. So it is enough to show that $N$ maps $N^{2p_r-2}v'_r$
($r\in[1,d]$) into $W_d$. But $NN^{2p_r-2}v'_r=N^{2p_r-1}v'_r=N^kN^{2p_r-1-k}v'_r\sub N^kV\sub W_d$. Thus 
$NW_d\sub W_d$. Hence $g^{*2}W_d=W_d$ and $\check{g^{*2}}W_d^\pe=W_d^\pe$. For $r\in[d+1,\s]$ we have 
$g\i v_r\in W_d^\pe$ by 1.7(iv). Since $\check{g^{*2}}W_d^\pe=W_d^\pe$ we have $g^{*j}v_r\in W_d^\pe$ for all odd
$j$; hence $gX_r\sub W_d^\pe$. Thus, $gW'_d\sub W_d^\pe$. Since $\dim W'_d=n-\dim W_d=\dim W_d^\pe$, it follows 
that $gW'_d=W_d^\pe$. Since $\check{g^{*2}}W_d^\pe=W_d^\pe$ it follows that $g^{*2}W'_d=W'_d$. Let $\d$ be the 
number of Jordan blocks of $N:W_d@>>>W_d$ that is, $\d=\dim(\ker N:W_d@>>>W_d)$. We have 
$\dim W_d-\d=\dim NW_d\ge\sum_{r\in[1,d]}(2p_r-2)=\dim W_d-d$. (The inequality follows from \cite{\WE, 3.1(b)}
applied to $N:W_d@>>>W_d$.) Hence $\d\le d$. From the definition of $\d$ we see that 
$\dim(\ker N^k:W_d@>>>W_d)\le\d k$. Recall that $\dim N^kW_d=\sum_{r\in[1,d]}(2p_r-1-k)=\dim W_d-kd$. Hence 
$\dim(\ker N^k:W_d@>>>W_d)=\dim W_d-\dim N^kW_d=kd$. Hence $kd\le\d k$. Since $k>0$ we deduce $d\le\d$. Hence 
$d=\d$. Since $\dim(\ker N^k:W_d@>>>W_d)=kd$ we see that each of the $d$ Jordan blocks of $N:W_d@>>>W_d$ has size
$\ge k$. Since $N^kW=N^kV$ and $V=W_d\op W'_d$ we see that $N^kW'_d=0$. Hence each Jordan block of
$N:W'_d@>>>W'_d$ has size $\le k$. This proves (b).

\subhead 5.3\endsubhead
In this subsection we assume that $g$ is as in 5.1. By 5.1(a) we can find $(V_*,V'_*)\in X_g$ so that the 
definitions and results in 4.2 are applicable. We show:

(a) {\it if $t\in[1,\s]$ then $((g^{*2}-1)^{p_t-1}x,g(g^{*2}-1)^{p_t-1}x)\ne0$ for some \lb
$x\in\ker(g^{*2}-1)^{2p_t-1}$.}
\nl
With the notation of 4.1 we have $(g^{*2})^k(z^t_0)=z^t_k$ if $k\in[0,2p_t-2]$,
$$(g^{*2})^{2p_t-1}(z^t_0)=\sum_{e\in[0,2p_t-2]}(-1)^e\bin{2p_t-1}{e}z^t_e$$
hence $(g^{*2}-1)^{2p_t-1}(z^t_0)=0$ and
$$(g^{*2}-1)^{p_t-1}(z^t_0)=\sum_{e\in[0,p_t-1]}(-1)^{p_t-1-e}\bin{p_t-1}{e}z^t_e.$$
Thus 
$$\align&((g^{*2}-1)^{p_t-1}z^t_0,g(g^{*2}-1)^{p_t-1}z^t_0)\\&=(\sum_{e\in[0,p_t-1]}(-1)^{p_t-1-e}\bin{p_t-1}{e}
z{}^t_e,\sum_{e\in[0,p_t-1]}(-1)^{p_t-1-e}\bin{p_t-1}{e}z'{}^t_e)\\&=(-1)^{p_t-1}(z{}^t_0,z'{}^t_{p_t-1})
=(-1)^{p_t-1}.\endalign$$
This proves (a).

\subhead 5.4\endsubhead
Let $\ct$ be the set of sequences $c_*=(c_1\ge c_2\ge c_3\ge\do)$ in $\NN$ such that $c_m=0$ for $m\gg0$ and 
$c_1+c_2+\do=\nn$. For $c_*\in\ct$ we define $c_*^*=(c^*_1\ge c^*_2\ge c^*_3\ge\do)\in\ct$ by 
$c^*_i=|\{j\ge1;c_j\ge i\}|$ and we set $\mu_i(c_*)=|\{j\ge1;c_j=i\}|$ ($i\ge1$); thus we have
$\mu_i(c_*)=c^*_i-c^*_{i+1}.$ For $i,j\ge1$ we have
$$i\le c_j\text{ iff }j\le c^*_i.\tag a$$
For $c_*,c'_*\in\ct$ we say that $c_*\le c'_*$ if the following (equivalent) conditions are satisfied:

(i) $\sum_{j\in[1,i]}c_j\le\sum_{j\in[1,i]}c'_j$ for any $i\ge1$;

(ii) $\sum_{j\in[1,i]}c^*_j\ge\sum_{j\in[1,i]}c'{}^*_j$ for any $i\ge1$.
\nl
The following result is proved in \cite{\LX, 1.4(e)}

(b) {\it Let $c_*,c'_*\in\ct$ and $i\ge1$ be such that $c_*\le c'_*$,
$\sum_{j\in[1,i]}c^*_j=\sum_{j\in[1,i]}c'{}^*_j$. Then $c^*_i\le c'{}^*_i$. If in addition $\mu_i(c_*)>0$, then
$\mu_i(c'_*)>0$.}

\subhead 5.5\endsubhead 
Until the end of 5.10 we assume that $p=2$ so that $G,D$ are as in 1.1(a),(b); then $g\in G^1$ is almost 
unipotent if and only if it is unipotent.
Let $u\in G^1$ be unipotent. We associate to $u$ the sequence $c_*\in\ct$ whose nonzero terms are the sizes of 
the Jordan blocks of $u^{*2}:V@>>>V$. We must have $\mu_i(c_*)=$even for any even $i$. We also associate to $u$ a
map $\e_u:\{i\in 2\NN+1;i\ne0,\mu_i(c_*)>0\}@>>>\{0,1\}$ as follows: $\e_u(i)=0$ if 
$((u^{*2}-1)^{(i-1)/2}x,u(u^{*2}-1)^{(i-1)/2}x)=0$ for all $x\in\ker(u^{*2}-1)^i:V@>>>V$ and $\e_u(i)=1$ 
otherwise; we have automatically $\e_u(i)=1$ if $\mu_i(c_*)$ is odd. Now $u\m(c_*,\e_u)$ defines a bijection 
$\uD@>\si>>\fS$ where $\fS$ is the set consisting of all pairs $(c_*,\e)$ where $c_*\in\ct$ is such that
$\mu_i(c_*)=$even for any even $i$ and $\e:\{i\in2\NN+1;i\ne0,\mu_i(c_*)>0\}@>>>\{0,1\}$ is a function such that 
$\e(i)=1$ if $\mu_i(c_*)$ is odd. (See \cite{\SP, I,2.7}). We denote by $\g_{c_*,\e}$ the element of $\uD$
corresponding to $(c_*,\e)\in\fS$. For $(c_*,\e)\in\fS$ it will be convenient to extend $\e$ to a function 
$\ZZ_{>0}@>>>\{-1,0,1\}$ (denoted again by $\e$) by setting $\e(i)=-1$ if $i$ is even or $\mu_i(c_*)=0$.

Now let $\g=\g_{c_*,\e},\g'=\g_{c'_*,\e'}$ with $(c_*,\e),(c'_*,\e')\in\fS$. Assume that 
$c_*=(2p_1-1,2p_2-1,\do,2p_\s-1,0,0,\do)$ and that $\e$ is such that $\e(2p_i-1)=1$ for $i\in[1,\s]$. Assume that
for some/any $g\in\g'$ we have $X_g\ne\em$. We will show that 

(a) {\it $\g$ is contained in the closure of $\g'$.}
\nl
The proof of (a) (given in 5.6-5.9) is almost a copy of the proof of Theorem 1.3 in \cite{\LX}. 

\subhead 5.6\endsubhead
It is enough to show that 

(a) $c_*\le c'_*$
\nl
and that for any $i\ge1$, (b),(c) below hold:

(b) $\sum_{j\in[1,i]}c^*_j-\max(\e(i),0)\ge\sum_{j\in[1,i]}c'{}^*_j-\max(\e'(i),0)$;

(c) if $\sum_{j\in[1,i]}c^*_j=\sum_{j\in[1,i]}c'{}^*_j$ and $c^*_{i+1}-c'{}^*_{i+1}$ is odd then $\e'(i)\ne0$.
\nl
(See \cite{\SP, II,8.2}.)

Let $g\in\g'$. We choose $(V_*,V'_*)\in X_g$. Then the notation and results in 5.2 are valid for $g$.
From 5.2(a) we see that (a) holds. Note also that, by \cite{\LX, 1.4(d)}, for $i\ge1$,

(d) {\it we have $\sum_{j\in[1,i]}c^*_j=\sum_{j\in[1,i]}c'{}^*_j$ iff $\dim N^iV=\L'_i$.}

\subhead 5.7\endsubhead
Let $i\ge1$. We show (compare \cite{\LX, 1.6}):

(a) {\it If $\mu_i(c_*)>0$ and $\sum_{j\in[1,i]}c^*_j=\sum_{j\in[1,i]}c'{}^*_j$ then $\e'(i)=1$.}
\nl
By 5.4(b) we have $\mu_i(c'_*)>0$. Since $\mu_i(c_*)>0$ we see that $i=2p_d-1$ for some $d\in[1,\s]$.
If $\mu_i(c'_*)$ is odd then $\e'(i)=1$ (by definition, since $i$ is odd). Thus we may assume that 
$\mu_i(c'_*)\in\{2,4,6,\do\}$. From our assumption we have that $\dim N^iV=\L'_i$ (see 5.6(d)).

Let $v_1,v_2,\do,v_\s$ be vectors in $V$ attached to $V_*,V'_*,g$ as in 4.7. For $r\in[1,\s]$ let $W_r,W'_r$ be 
as in 5.2; we set $W_0=0$, $W'_0=V$. From 5.2(b) we see that $N^iW'_{d-1}=0$ at least if $d\ge2$; but the same
clearly holds if $d=1$. We have $g^{*(-2p_d)}v_d\in W'_{d-1}$; since $W'_{d-1}$ is $g^{*2}$-stable (see 5.2) we
have $v_d\in W'_{d-1}$ hence $N^{2p_d-1}v_d=0$ (see 5.2) and
$$\align&(N^{p_d-1}v_d,gN^{p_d-1}v_d)\\&=\sum_{e,e'\in[0,p_d-1]}\bin{p_d-1}{e}\bin{p_d-1}{e'}
(g^{*(2e)}v_d,g^{*(2e'+1)}v_d)\\&=(v_d,g^{*(2p_d-1)}v_d)=1.\endalign$$
(We have used that $(v_d,g^{*j}v_d)=0$ if $j\in[-2p_d+2,2p_d-2]_{odd}$ and \lb
$(v_d,g^{*(2p_d-1)}v_d)=1$.)
Thus $\e'(i)=1$. This proves (a).

\subhead 5.8\endsubhead
We prove 5.6(b). It is enough to show that, if $\e(i)=1$ and $\e'(i)\le0$ then
$\sum_{j\in[1,i]}c^*_j\ge\sum_{j\in[1,i]}c'{}^*_j+1$. Assume this is not so. Then using 5.6(a) we have
$\sum_{j\in[1,i]}c^*_j=\sum_{j\in[1,i]}c'{}^*_j$. Since $\e(i)=1$ we have $\mu_i(c_*)>0$; using 5.7(a) we see 
that $\e'(i)=1$, a contradiction. Thus 5.6(b) holds.

\subhead 5.9\endsubhead
We prove 5.6(c). If $i$ is even then $\e'(i)=-1$, as required. Thus we may assume that $i$ is odd. Using 5.6(a) 
and 5.4(b) we see that $c^*_i\le c'{}^*_i$. Assume first that $c^*_i=c'{}^*_i$. From 
$\mu_i(c_*)=c^*_i-c^*_{i+1}$, $\mu_i(c'_*)=c'{}^*_i-c'{}^*_{i+1}$ we deduce that 
$\mu_i(c_*)-\mu_i(c'_*)=c'{}^*_{i+1}-c^*_{i+1}$ is odd. If $\mu_i(c'_*)$ is odd we have $\e'(i)=1$ (since $i$ is 
even); thus we have $\e'(i)\ne0$, as required. If $\mu_i(c'_*)=0$ we have $\e'(i)=-1$; thus we have $\e'(i)\ne0$,
as required. If $\mu_i(c'_*)\in\{2,4,6,\do\}$ then $\mu_i(c_*)$ is odd so that $\mu_i(c_*)>0$ and then 5.7(a) 
shows that $\e'(i)=1$; thus we have $\e'(i)\ne0$, as required. 

Assume next that $c^*_i<c'{}^*_i$. By 5.6(a) we have $\sum_{j\in[1,i+1]}c^*_j\ge\sum_{j\in[1,i+1]}c'{}^*_j$; 
using the assumption of 5.6(c) we deduce that $c^*_{i+1}\ge c'{}^*_{i+1}$. Combining this with $c^*_i<c'{}^*_i$ 
we deduce $c^*_i-c^*_{i+1}<c'{}^*_i-c'{}^*_{i+1}$ that is, $\mu_i(c_*)<\mu_i(c'_*)$. It follows that 
$\mu_i(c'_*)>0$. If $\mu_i(c_*)>0$ then by 5.7(a) we have $\e'(i)=1$; thus we have $\e'(i)\ne0$, as required. 
Thus we can assume that $\mu_i(c_*)=0$. We then have $c^*_i=c^*_{i+1}$ and we set $\d=c^*_i=c^*_{i+1}$. As we 
have seen earlier, we have $c^*_{i+1}\ge c'{}^*_{i+1}$; using this and the assumption of 4.6(c) we see that 
$c^*_{i+1}-c'{}^*_{i+1}=2a+1$ where $a\in\NN$. It follows that $c'{}^*_{i+1}=\d-(2a+1)$. In particular we have
$\d\ge2a+1>0$.

If $k\in[0,2a]$ we have $c'_{\d-k}=i$. (Indeed, assume that $i+1\le c'_{\d-k}$; then by 5.4(a) we have 
$\d-k\le c'{}^*_{i+1}=\d-(2a+1)$ hence $k\ge 2a+1$, a contradiction. Thus $c'_{\d-k}\le i$. On the other hand,
$\d=c^*_i<c'{}^*_i$ implies by 5.4(a) that $i\le c'_\d$. Thus $c'_{\d-k}\le i\le c'_\d\le c'_{\d-k}$ hence
$c'_{\d-k}=i$.)

Using 5.4(a) and $c'{}^*_{i+1}=\d-(2a+1)$ we see that $c'_{\d-(2a+1)}\ge i+1$ (assuming that $\d-(2a+1)>0$). Thus
the sequence $c'_1,c'_2,\do,c'_\d$ contains exactly $2a+1$ terms equal to $i$, namely
$c'_{\d-2a},\do,c'_{\d-1},c'_\d$.

We have $i>c_{\d+1}$. (If $i\le c_{\d+1}$ then from 5.4(a) we would get $\d+1\le c^*_i=\d$, a contradiction.)

Since $\d>0$, from $c^*_i=\d$ we deduce that $i\le c_\d$ (see 5.4(a)); since 
$\mu_i(c_*)=0$ we have $c_\d\ne i$ hence $c_\d>i$. From the assumption of 5.6(c) we see that $\dim N^iV=\L'_i$ 
(see 5.6(d)). Using this and $c_\d>i>c_{\d+1}$ we see that 5.2(b) is applicable and gives that 
$V=W_\d\op W'_\d$, $W_\d,W'_\d$ are $g^{*2}$-stable and $gW'_\d=W_\d^\pe$; moreover, $g^{*2}:W_\d@>>>W_\d$ has 
exactly $\d$ Jordan blocks and each one has size $\ge i$ and $g^{*2}:W'_\d@>>>W'_\d$ has only Jordan blocks of 
size $\le i$. Since the $\d$ largest numbers in the sequence $c'_1,c'_2,\do$ are $c'_1,c'_2,\do,c'_\d$ we see 
that the sizes of the Jordan blocks of $g^{*2}:W_\d@>>>W_\d$ are $c'_1,c'_2,\do,c'_\d$. Since the last sequence 
contains an odd number ($=2a+1$) of terms equal to $i$ we see that $\e_{g|_{W_\d}}(i)=1$. (We explain the meaning
of the last equation. The dual space $W_\d^*$ of $W_\d$ can be identified canonically with $V^*/W_\d^\pe$ which 
can be identified canonically with $gW_\d$ since $V^*=gW_\d\op gW'_\d=gW_\d\op W_\d^*$. Hence $x\m gx$ can be 
viewed as an isomorphism $W_\d@>>>W_\d^*$, denoted by $g|_{W_\d}$. Note that $(g|_{W_\d})^{*2}$ is equal to the 
restriction of $g^{*2}$ to $W_\d$ hence is unipotent hence $g|_{W_\d}$ is unipotent and its invariant 
$\e_{g|_{W_\d}}(i)$ is defined.) Hence there exists $z\in W_\d$ such that $(g^{*2}-1)^iz=0$ and 
$((g^{*2}-1)^{(i-1)/2}z,g(g^{*2}-1)^{(i-1)/2}z)=1$. This shows that $\e_g(i)=1$ that is $\e'(i)=1$. This 
completes the proof of 5.6(c) and also completes the proof of 5.5(a).   

From 5.5(a) (which is applicable in view of 5.3(a)) we see that 1.3(a) holds for any $C\in\uWW^{wl}$. Moreover we
see that 1.3(b), 1.3(d) hold for $G,D$ (we use the description of distinguished unipotent classes of $D$ given in
\cite{\SP}). This completes the verification of 1.4 for our $G,D$.

\subhead 5.10\endsubhead
We now prove Theorem 1.14 for our $G,D$. (Recall that $p=2$). We can assume that $w=z_{p_*}$ with
$p_*=(p_1\ge p_2\ge\do\ge p_\s)$ as in 4.4; thus $2(p_1+p_2+\do+p_\s)=n+\s$. We then have 
$$\ul(w)=p_1+3p_2+\do+(2\s-1)p_\s-(\s^2-\s)/2.$$
Since $g^{*2}$ has Jordan blocks of sizes $2p_1-1,2p_2-1,\do,2p_\s-1$, we see from \cite{\SP, p.96} and 5.3(a)
that
$$d:=\dim(Z(g))=\sum_{h\ge1}f_h^2/2-\sum_{h\ge1,\text{ odd}}f_h+n/2$$
where $f_1\ge f_2\ge\do$ is the partition dual to $2p_1-1,\do,2p_\s-1$; thus, for $j\ge1$ we have
$f_j=\sha(i\in[1,\s];2p_i-1\ge j)$. Note also that $\cz_{G^0}^D=\{1\}$. We must show that $\ul(w)=d$.
We can find integers $a_1,a_2,\do,a_t,b_1,b_2,\do,b_t$ (all $\ge1$) such that $b_1+b_2+\do+b_t=\s$,
$2p_i-1=a_1+a_2+\do+a_t$ for $i\in[1,b_1]$, $2p_i-1=a_1+a_2+\do+a_{t-1}$ for $i\in[b_1+1,b_1+b_2]$, $\do$, 
$2p_i-1=a_1$ for $i\in[b_1+b_2+\do+b_{t-1}+1,b_1+b_2+\do+b_t]$. We have 
$$\align&2\ul(w)=(a_1+\do+a_t+1)b_1^2+(a_1+\do+a_{t-1}+1)((b_1+b_2)^2-b_1^2)+\do\\&+(a_1+1)((b_1+\do+b_t)^2
-(b_1+\do+b_{t-1})^2)-(\s^2+\s)\\&=a_1(b_1+\do+b_t)^2+a_2(b_1+\do+b_{t-1})^2+\do+a_tb_1^2-\s.\endalign$$
We have 
$$\align&2d=a_1(b_1+\do+b_t)^2+a_2(b_1+\do+b_{t-1})^2+\do+a_tb_1^2\\&
-(a_1+1)(b_1+\do+b_t)-a_2(b_1+\do+b_{t-1})-\do-a_tb_1+2(n/2).\endalign$$
To show that $d=\ul(w)$, we must show:
$$\align&a_1(b_1+\do+b_t)^2+a_2(b_1+\do+b_{t-1})^2+\do+a_tb_1^2-\s\\&=a_1(b_1+\do+b_t)^2+a_2(b_1+\do+b_{t-1})^2
+\do+a_tb_1^2\\&-(a_1+1)(b_1+\do+b_t)-a_2(b_1+\do+b_{t-1})-\do-a_tb_1+n.\endalign$$
or that
$$n=a_1(b_1+\do+b_t)+a_2(b_1+\do+b_{t-1})+\do+a_tb_1.$$
But this is just another way to write the identity $n=(2p_1-1)+(2p_2-1)+\do+(2p_\s-1)$. This completes the proof.

\subhead 5.11\endsubhead
In this subsection we assume that $p\ne2$. Let $D_*=\{s\in D;s^{*2}=1\}$. Now $D_*$ is in bijection with the set 
of nondegenerate symmetric bilinear forms on $V$ with values in $\kk$: to $s\in D_*$ corresponds the symmetric 
bilinear form $\us$ given by $x,x'\m s(x)(x')=s(x')(x)$. Note that for $g\in D$ we have $g\in D_{au}$ if and 
only if $g=su$ where $s\in D_*$ and $u$ is a unipotent element in $G^0$ commuting with $s$ or, equivalently, $u$ 
is a unipotent element of the orthogonal group $O_{\us}$. It follows that the set $\un{D_{au}}$ of 
$G^0$-conjugacy classes in $D_{au}$ is in bijection with the set of partitions of $n$ such that each even part 
appears an even number of times: to $\g\in\un{D_{au}}$ corresponds the partition given by the sizes of the Jordan
blocks of $g^{*2}:V@>>>V$ for some/any $g\in\g$. We will define a canonical map $\Ph:\uWW@>>>\un{D_{au}}$ using 
the same principle as that used in Theorem 1.3. For $w\in\WW,\g\in\un{D_{au}}$ we write $w\dsv_D\g$ if for 
some/any $g\in\g$ and some $B\in\cb$ we have $(B,gBg\i)\in\co_w$. For $w\in\WW$ we set 
$\Si'_{w,D}=\{\g\in\un{D_{au}};w\dsv_D\g\}$; we regard $\Si'_{w,D}$ as a partially ordered set where $\g\le\g'$ if
$\g\sub\bar\g'$ (closure of $\g'$ in $D$). Using arguments as in the proof of 1.2(a) we see that if $w,w'$ are 
elements of $\WW_{D-min}$ which are $\e_D$-conjugate, then $\Si'_{w,D}=\Si'_{w',D}$. Hence for any $C\in\uWW_D$ 
we can define $\Si'_{C,D}=\Si'_{w,D}$ where $w$ is any element of $C_{min}$. We have the following result.

(a) {\it Let $C\in\uWW_D^{el}$. There exists (a necesarily unique) $\g\in\Si'_{C,D}$ such that $\g\le\g'$ for all 
$\g'\in\Si'_{C,D}$. We set $\g=\Ph'(C)$.}
\nl
We can assume that $C$ contains $z_{p_*}$ (as in 4.4) so that $\Si'_{C,D}=\Si'_{z_{p_*},D}$. Then the result 
follows from 5.1 and 5.2(a). Note that for $C$ as above, $\Ph'(C)$ is the $G^0$-orbit in $D_{au}$ corresponding to
the partition $2p_1-1,2p_2-1,\do,2p_\s-1$ of $n$. We note also that the analogues of 1.3(b),(d) clearly hold in 
our case. Namely, if $C,C'$ are elements of $\uWW_D^{el}$ and $\Ph'(C)=\Ph'(C')$, then $C=C'$; moreover if
$\g\in\un{D_{au}}$ is distinguished (in the sense that the centralizer in $G^0$ of an element of $\g$ contains no
torus $\ne1$ then $\g=\Ph'(C)$ for some $C\in\uWW_D^{el}$. Also if $C\in\uWW_D^{el}$, $w\in C_{min}$ and 
$\g=\Ph'(C)$ then for any $g\in\g$ we have $\dim(Z_G(g))=\ul(w)$ (compare 1.14); the proof is the same as that in
5.10. Using arguments similar to those in the proof of 1.4(a) we see that (a) holds in the same form for any 
$C\in\uWW_D$ (not necessarily elliptic). Thus, $\Ph':\uWW_D@>>>\un{D_{au}}$ is well defined and surjective. 

We now give a combinatorial description of the map $\Ph':\uWW_D>>>\un{D_{au}}$. We identify $\uWW_D$ with the set
of partitions of $n$ (to the $\e_D$-conjugacy class of $w\in\WW$ corresponds the partition which gives the 
lengths of the cycles of the permutation $w\ww$ of $[1,n]$, $\ww$ as in 4.3). We identify $\un{D_{au}}$ with the
set of partitions of $n$ in which any even part appears an even number of times (to the $GL(V)$-conjugacy class 
of $g\in D_{au}$ we associate the partition which gives the sizes of the Jordan blocks of $g^{*2}$). Then
$\Ph'(\l)=\z$ where each odd part $2a+1$ of $\l$ gives a part of size $2a+1$ of $\z$ and each even part $2a$ of
$\l$ gives two parts of size $a,a$ of $\z$. For example if $\l=(5,4,3,3,2,2,1,1)$ then
$\Ph'(\l)=(5,3,3,2,2,1,1,1,1,1,1)$.

From this combinatorial description we see that we have the following analogue of Theorem 1.16.

{\it For any $\g\in\un{D_{au}}$ the function $\Ph'{}\i(\g)@>>>\NN$, $C\m\mu(C)$ ($\mu$ defined as in the proof 
of 1.4(a)) reaches its minimum at a unique element $C_0\in\Ph'{}\i(\g)$. Thus we have a well defined map 
$\Ps':\un{D_{au}}@>>>\uWW_D$, $\g\m C_0$ such that $\Ph'\Ps':\un{D_{au}}@>>>\un{D_{au}}$ is the identity map.}

\subhead 5.12\endsubhead
In this subsection $p$ is arbitrary. We fix a sequence $p_*=(p_1\ge p_2\ge\do\ge p_\s)$ in $\ZZ_{>0}$ such that 
$n=(2p_1-1)+\do+(2p_\s-1)$ and $(V_* ,V'_*)\in\co_w$ where $w=z_{p_*}$. Let $g\in G^1$ be such that 
$g\cdot V_*=V'_*$ and such that $N:=g^{*2}-1:V@>>>V$ is nilpotent with Jordan blocks of sizes 
$2p_1-1,2p_2-1,\do,2p_\s-1$. Let $v_1,v_2,\do,v_\s$ be the sequence of vectors associated to $g$ in 4.7 (with 
$u=\s$); these vectors are defined up to multiplication by $\pm1$. For $t\in[1,\s]$ and $i\in\ZZ$ we set 
$z^t_i=g^{*(2i-2p_t)}v_t$ and $z'{}^t_i=gz^t_i$. From the definitions we see that 

(a) $z'{}^t_i=gz^t_i,z^t_{i+1}=\chg z'{}^t_i$ for all $t,i$;

(b) $(z^t_i,z'{}^t_j)=0$ if $j-i\in[-p_t+1,p_t-2]$, $(z^t_i,z'{}^t_j)=1$ if $j-i\in\{-p_t,p_t-1\}$ ($t\in[1,\s]$,
$i,j\in\ZZ$).
\nl
Moreover, by the argument in 4.8(a) we have that $\{z^t_i;t\in[1,\s],i\in[0,2p_t-2]\}$ is a basis of $V$. For any 
$t\in[1,\s]$ let $X_t$ be the subspace of $V$ spanned by $\{z^t_i;i\in[0,2p_t-2]\}$ so that
$V=X_1\op X_2\op\do X_\s$. Note that for any $k\ge0$ we have $\dim N^k(V)=\sum_{r\in[1,\s]}\max(2p_r-1-k,0)$.
Applying 5.2(b) with $k=2p_d-1$ for $d=1,2,\do,\s$ we see that each of the subspaces 
$$X_1\sub X_1\op X_2\sub\do\sub X_1\op X_2\op\do\op X_\s,$$
$$X_2\op\do\op X_{\s+1}\sps\do\sps X_\s\op X_{\s+1}\sps X_{\s+1}$$ 
of $V$ is $g^{*2}$-stable. Taking intersections we see that each of the subspaces $X_1,X_2,\do,X_\s$ of $V$ is 
$g^{*2}$-stable. From 5.2(b) we see also that 
$$(X_1\op X_2\op\do\op X_e,g(X_{e+1}\op X_{e+2}\op\do\op X_\s))=0$$ 
for any $e\in[1,\s]$ hence $(X_t,gX_r)$ for $t<r$ in $[1,\s]$. Using this and 4.1(a) we see that
if $t<r$ then $(X_r,gX_t)=(X_t,g^{*(-1)}X_r)=(X_t,gX_r)=0$. (We have use that $X_r$ is $g^{*2}$-stable.) Thus
$(X_r,gX_t)=0$ for any $t\ne r$ in $[1,\s]$. We see that 

(c) $(z^r_i,z'{}^t_j)=0$ if $t\ne r$, $i,j\in\ZZ$.
\nl
We now show for any $r\in[1,\s]$ that

(d) $(z^r_i,z'{}^r_j)=\bin{j-i+p_r-1}{j-i-p_r+1}$ if $j-i\ge p_r-1$,   

(e) $(z^r_i,z'{}^r_j)=\bin{i-j+p_r-2}{i-j-p_r}$ if $j-i\le-p_r$.       
\nl
We first prove (d). When $j-i=p_r-1$, (d) holds by (b). Now assume that $j-i=p_r-1+s$ where $s>0$ and that (d) is
known when $j-i=p_r-1+s'$, $0\le s'<s$. Since $g^{*2}$ restricts to a unipotent automorphism of $X_r$ and 
$\dim X_r=2p_r-1$, we see that $N^{2p_r-1}$ acts as $0$ on $X_r$, so that $(g^{*2}-1)^{2p_r-1}z^r_0=0$. Since 
$g^{*2}z^r_0=z^r_1$, $g^{*2}z^r_1=z^r_2$, $\do$, $g^{*2}z^r_{2p_r-2}=z^t_{2p_r-1}$, it follows that 
$$\sum_{k\in[0,2p_r-1]}(-1)^k\bin{2p_r-1}{k}z^r_k=0.$$
$$\sum_{k\in[0,2p_r-1]}(-1)^k\bin{2p_r-1}{k}(z^r_k,z'{}^r_{p_r-1+s})=0.$$
Using (b) and the induction hypothesis we deduce
$$\sum_{k\in[1,2p_r-1]}(-1)^k\bin{2p_r-1}{k}\bin{p_r-1+s-k+p_r-1}{s-k}+(z^r_0,z'{}^r_{p_r-1+s})=0.$$
It is then enough to show that
$$\sum_{k\in[0,2p_r-1]}(-1)^k\bin{2p_r-1}{k}\bin{p_r-1+s-k+p_r-1}{s-k}=0$$
for $s>0$ or setting $m=s-k$ that
$$\sum_{s\ge0}\sum_{k\in[0,2p_r-1];m\ge0;k+m=s}(-1)^k\bin{2p_r-1}{k}\bin{2p_r+m-2}{m}T^s=1$$
where $T$ is an indeterminate. An equivalent statement is
$$(\sum_{k\in[0,2p_r-1]}(-1)^k\bin{2p_r-1}{k}T^k)(\sum_{m\ge0}\bin{2p_r+m-2}{m}T^m)=1.$$
This folows from the identity $\sum_{m\ge0}\bin{M+m-1}{m}T^m=(1-T)^{-M}$ (for $M\ge1$) which is easily verified.
This proves (d).

We now prove (e). Assume that $j-i\le-p_r$. Using 4.1(a), we have
$$(z^r_i,z'{}^r_j)=(z^r_i,gz^r_j)=(z^r_j,g^{*(-1)}z^r_i)=(z^r_j,gz^r_{i-1})=(z^r_j,z'{}^r_{i-1}).$$
Note that $i-1-j\ge p_r-1$ hence using (d) we have 
$$(z^r_j,z'{}^r_{i-1})=\bin{i-1-j+p_r-1}{i-1-j-p_r+1}$$
and (e) follows.

Now let $\tg\in G^1$ be another element such that $\tg\cdot V_*=V'_*$ and such that $\ti N:=\tg^{*2}-1:V@>>>V$ is
nilpotent with Jordan blocks of sizes $2p_1-1,2p_2-1,\do,2p_\s-1$. We associate to $\tg$ vectors 
$\tz^t_i\in V,\tz'{}^t_i\in V^*$ in the same way as $z^t_i\in V,z'{}^t_i\in V^*$ were associated to $g$. From 
(a),(b),(c),(d),(e) for $g$ and $\tg$ we see that

$z'{}^t_i=gz^t_i,z^t_{i+1}=\chg z'{}^t_i$, $\tz'{}^t_i=\tg\tz^t_i,\tz^t_{i+1}=\che{\tg}\tz'{}^t_i$ for all $t,i$;

$(z^t_i,z'{}^t_j)=(\tz^t_i,\tz'{}^t_j)=0$ if $j-i\in[-p_t+1,p_t-2]$;

$(z^t_i,z'{}^t_j)=(\tz^t_i,\tz'{}^t_j)=\bin{j-i+p_t-1}{j-i-p_t+1}$ if $j-i\ge p_t-1$;   

$(z^t_i,z'{}^t_j)=(\tz^t_i,\tz'{}^t_j)=\bin{i-j+p_t-2}{i-j-p_t}$ if $j-i\le-p_t$;       

$(z^r_i,z'{}^t_j)=(\tz^r_i,\tz'{}^t_j)=0$ if $t\ne r$, $i,j\in\ZZ$.
\nl
Since $\{z^t_i;t\in[1,\s],i\in[0,2p_t-2]\}$ and $\{\tz^t_i;t\in[1,\s],i\in[0,2p_t-2]\}$ are bases of $V$ there is
a unique isomorphism of vector spaces $T:V@>>>V$ such that $\tz^t_i=T(z^t_i)$ for all 
$t\in[1,\s],i\in[0,2p_t-1]$. The formulas above show that $(z^r_i,gz^t_j)=(T(z^r_i),\tg T(z^t_j))$ (that is, 
$(z^r_i,gz^t_j)=(z^r_i,\che{T}\i\tg T(z^t_j))$) for any $r,t$ in $[1,\s]$ and any 
$i\in[0,2p_t-2],j\in[0,2p_r-2]$. It follows that $g=\che{T}\i\tg T$.

From the definitions, for any $r\in[1,\s],i\in[1,p_r]$, the subspace $V_{p_1+\do+p_{r-1}+i}$ is generated by 
$z^t_h$ ($t<r$, $h\in[0,p_h-1]$) and by $z^r_h(h\in[0,i-1])$; similarly, the subspace $V_{p_1+\do+p_{r-1}+i}$ is 
generated by $\tz^t_h$ ($t<r$, $h\in[0,p_h-1]$) and by $\tz^r_h (h\in[0,i-1])$. Also, for any $r\in[1,\s]$ and
any $i\in[1,p_r-1]$, the subspace $V'_{(p_1-1)+\do+(p_{r-1}-1)+i}$ is generated by $z^t_h$ ($t<r,h\in[1,p_h-1]$) 
and by $z^r_h$ ($h\in[1,i])$; similarly, the subspace $V'_{(p_1-1)+\do+(p_{r-1}-1)+i}$ is generated by $\tz^t_h$ 
($t<r$, $h\in[1,p_h-1]$) and by $\tz^r_h(h\in[1,i])$. It follows that $TV_u=V_u$ for $u\in[0,p_1+\do+p_\s]$ and 
$TV'_u=V'_u$ for $u\in[1,p_1+\do+p_\s-\s]$. For any $i\in[0,n]$ we have 
$$V'_i=(gV_{n-i})^\pe=(\tg V_{n-i})^\pe=(\che T gT\i V_{n-i})^\pe=T((gT\i V_{n-i})^\pe).$$ 
For $u\in[p_1+\do+p_\s-\s,n-1]$ we have $n-u\in[1,p_1+\do+p_\s]$ hence $T\i V_{n-u}=V_{n-u}$ and
$V'_u=T((gV_{n-u})^\pe)=TV'_u$. Thus $TV'_i=V'_i$ for all $i\in[0,n]$. For $u\in[p_1+\do+p_\s,n-1]$ we have
$$gV_u=\tg V_u=\che(T)gT\i V_u=(V'_{n-u})^\pe=(TV'_{n-u})^\pe=\che{T}((V'_{n-u})^\pe)=\che{T}gV_u$$
hence $T\i V_u=V_u$. Thus $TV_i=V_i$ for all $i\in[0,n]$. Thus we have the following result (where $\g$ denotes
the set of all $g\in D_{au}$ such that the Jordan blocks of $g^{*2}-1:V@>>>V$ have sizes 
$2p_1-1,2p_2-1,\do,2p_\s-1$): 

(f) {\it The $GL(V)$-action $x:(g,B)\m(xgx\i,xBx\i)$ on $\{(g,B)\in\g\T\cb;(B,gBg\i)\in\co_w\}$ is transitive.}
\nl
In the special case where $p=2$ this proves Theorem 1.15 for our $G,D$.

\head 6. Finite unitary groups\endhead
\subhead 6.1\endsubhead
In this section we preserve the setup of 4.1 and we assume that $q>1$. We fix $g\in G^1$. Then 
$\ph:=g^{*2}=\chg g:V@>>>V$ (resp. $\ph':=g\chg:V^*@>>>V^*$) is the Frobenius map for a rational structure on $V$
(resp. $V^*$) over the finite subfield $F_{q^2}$ with $q^2$ elements of $\kk$. Note that 
$V^\ph:=\{x\in V;\ph(x)=x\}$, $V^{*\ph'}:=\{\x\in V^*;\ph'(\x)=\x\}$ are $F_{q^2}$-vector spaces of dimension $n$
such that $V=\kk\ot_{F_{q^2}}V^\ph$, $V^*=\kk\ot_{F_{q^2}}V^{*\ph'}$ and that $g$ restricts to an $F_q$-linear 
isomorphism $V^\ph@>>>V^{*\ph'}$. For $x,x'\in V$ we set $\la x,x'\ra=(x,gx')$. From 4.1(a) we deduce
$\la \ph(x),x'\ra=\la x',x\ra^q$ for any $x,x'\in\kk$. In particular, if $x,x'\in V^\ph$ we have
$$\la x,x'\ra=\la x',x\ra^q.\tag a$$
Applying (a) twice we see that $\la x,x'\ra=\la x,x'\ra^{q^2}$ for any $x,x'\in V^\ph$. Thus $x,x'\m\la x,x'\ra$ 
is a (nondegenerate) hermitian form $V^\ph\T V^\ph@>>>F_{q^2}$. It follows that we can find a basis
$e_1,e_2,\do,e_n$ of $V$ such that $\ph(e_i)=e_i$ for $i\in[1,n]$ and $(e_i,ge_j)=\d_{ij}$ for $i,j\in[1,n]$.

We define a group isomorphism $\ps:G^0@>>>G^0$ by $x\m g*x*g^{*(-1)}=\chg\check{x}\chg\i$. We express $\ps$ in
coordinates. Let $e'_1,e'_2,\do,e'_n$ be the basis of $V^*$ such that $(e_i,e'_j)=\d_{ij}$ for $i,j\in[1,n]$. 
Note that $e'_i=ge_i$ for all $i$. Let $x\in G^0$. We have $xe_i=\sum_jx_{ij}e_j$, $x_{ij}\in\kk$. Let 
$(x'_{ij})$ be the matrix which is transpose inverse to $(x_{ij})$. Then $\ps(x)(e_i)=\sum_jx'_{ij}{}^qe_j$ for 
all $i$. We see that $\ps$ is the Frobenius map for an $F_q$-rational structure on $G^0$. As pointed out in 4.2 
for any $V_*\in\cf$ we have $B_{g\cdot V_*}=g*B_{V_*}*g^{*(-1)}=\ps(B_{V_*})$. Hence, if $p_*,w=z_{p_*}$ are as 
in 4.4, the condition that $(V_*,g\cdot V_*)\in\co_w$  is equivalent to the condition that 
$B_{V_*}\in\XX_{z_{p_*}}$ where $\XX_{z_{p_*}}$ is the variety of all $B\in\cb$ such that $B,\ps(B)$ are in 
relative position $w$ (see \cite{\DL}). Thus $(V_*,g\cdot V_*)\m B_{V_*}$ is an isomorphism 
$X_g@>\si>>\XX_{z_{p_*}}$. Composing with 4.9(a) we deduce that $\XX_{z_{p_*}}$ is isomorphic to the variety 
consisting of all sequences $L_1,L_2,\do,L_\s$ of lines in $V$ such that for any $r\in[1,\s]$ we have

$\la L_r,\ph^hL_t\ra=0$ if $t\in[1,\s]$, $t<r$, $h\in[-p_t,p_t-2]$;

$\la L_r,\ph^hL_r\ra=0$ if $h\in[-p_r+1,p_r-2]$;

$\la L_r,\ph^{p_r-1}L_r\ra\ne0$.
\nl
Thus $\XX_{z_{p_*}}$ can be viewed as the quotient of the variety consisting of all sequences
$v_1,v_2,\do,v_\s$ of vectors in $V$ such that for any $r\in[1,\s]$ we have

$\la v_r,\ph^hv_t\ra=0$ if $t\in[1,\s]$, $t<r$, $h\in[-p_t,p_t-2]$;

$\la v_r,\ph^hv_r\ra=0$ if $h\in[-p_r+1,p_r-2]$;

$\la v_r,\ph^{p_r-1}v_r\ra=1$,
\nl
(a model of the variety $\ti\XX_{z_{p_*}}$ of \cite{\DL}) by the obvious (free) action of the group 
$$\{(\l_1,\do,\l_\s)\in\kk^{*\s};\l_r^{q^{2p_r-1}+1}=1 \text{ for all }r\}.$$
Hence this description provides a new proof of the affineness of $\XX_{z_{p_*}}$ and $\ti\XX_{z_{p_*}}$ 
(the first proof is given in \cite{\HEA}; another proof is given in \cite{\HL}).

\head 7. Final remarks\endhead
\subhead 7.1\endsubhead
In this section we assume that $G$ is as in 1.1(a). Let $\bG_p$ be the set of unipotent elements of $\bG:=G/G^0$. 
For any $D\in\bG_p$ we will write $\Ph_D:\uWW_D@>>>\uD$ for the map denoted by $\Ph$ in 1.3(a) (relative to 
$D$). If $\a:G@>>>G$ is an automorphism of algebraic groups and $D\in\bG_p$ then $\a$ induce naturally bijections
$\a_*:\uWW_D@>>>\uWW_{\a(D)}$, $\a_*:\uD@>>>\un{\a(D)}$ and from the definitions we have
$$\a_*(\Ph_D((C))=\Ph_{\a(D)}(\a_*(C))\tag a$$ 
for any $C\in\uWW_D$. We consider the semidirect product $\bG\WW$ in which for $D\in\bG$ and $w\in\WW$ we have 
$DwD\i=\e_D(w)$. Let $c_p(\bG)$ be the set of conjugacy classes in $\bG\WW$ of elements of the form $Dw$, 
$D\in\bG_p,w\in\WW$. Let $\uG$ be the set of unipotent $G$-conjugacy classes in $G$. Let $D\in\bG_p,w\in\WW$. Let
$C$ be the $\e_D$-conjugacy class of $w$. Let $\g_{Dw}$ be the unipotent $G$-conjugacy class of $G$ that contains
$\Ph_D(C)$. We show that for $D'\in\bG$, $w'\in\WW$ we have
$$\g_{D'w'Dw(D'w')\i}=\g_{Dw}.\tag b$$
We have
$$\align&D'w'Dw(D'w')\i\\&=D'D\e_D\i(w')ww'{}\i D'{}\i=D'D D'{}\i\e_{D'}(\e_D\i(w'))\e_{D'}(w)\e_{D'}(w'{}\i)
\\&=D'DD'{}\i(\e_{D'}\e_D\i\e_{D'{}\i}(\e_{D'}(w')))\e_{D'}(w)\e_{D'}(w')\i).\endalign$$
Thus $D'w'Dw(D'w')\i$ is $\e_{D'DD'}$-conjugate to $\e_{D'}(w)$. Let $C'$ be the $\e_{D'DD'{}\i}$-conjugacy class
of $\e_{D'}(w)$. It is enough to show that $\Ph_D(C),\Ph_{D'DD'{}\i}(C')$ are contained in the same $G$-conjugacy
class of $G$. Let $\a:G@>>>G$ be the automorphism $x\m g'xg'{}\i$ where $x\in D'$ is fixed. We have 
$\a(D)=D'DD'{}\i$ and $\a_*(C)=C'$. Using (a) we have 
$$\Ph_{D'DD'{}\i}(C')=\Ph_{\a(D)}(\a_*(C))=\a_*(\Ph_D(C))=x\Ph_D(C)x\i.$$
Thus $\Ph_D(C),\Ph_{D'DD'{}\i}(C')$ are indeed contained in the same $G$-conjugacy class of $G$. This proves (b).

From (b) we see that $Dw\m\g_{Dw}$ is constant on the $\bG\WW$-conjugacy classes of $\bG\WW$ which comprise
$c_p(\bG\WW)$; hence it defines a map $\hat\Ph:c_p(\bG\WW)@>>>\uG$ which is surjective since $\Ph_D$ is 
surjective for any $D\in\bG_p$.

\enddocument